\documentclass[oneside,a4paper]{article}
\usepackage[a4paper, margin=1.3in]{geometry}

\usepackage{authblk}
\usepackage{graphicx}
\graphicspath{{./}}
\DeclareGraphicsExtensions{.png,.pdf,.jpg,.jpeg}
\usepackage{placeins} 
\usepackage{adjustbox} 
\usepackage{float} 
\usepackage{tabularx} 
\usepackage[nottoc
]{tocbibind}
\usepackage{hyperref}
\usepackage{amsmath}
\usepackage{amsthm}
\usepackage{amssymb}
\usepackage[capitalize]{cleveref}
\usepackage{autonum}
\usepackage{eufrak}
\usepackage{xcolor}
\usepackage[normalem]{ulem}
\usepackage{bbm}
\usepackage{tikz}

\usepackage[ruled,vlined]{algorithm2e}

\theoremstyle{plain}

\theoremstyle{remark}
\newtheorem{rmrk}{Remark}
\theoremstyle{definition}

\newcommand{\dd}{{\mathbf d}}
\newcommand{\dbar}{{\mathbf{g}}}
\newcommand{\dmapt}{{\mathbf T_t}}
\newcommand{\dmaptinv}{{\mathbf T_t^{-1}}}
\newcommand{\uu}{{\mathbf u}}
\newcommand{\vv}{{\mathbf v}}
\newcommand{\xx}{{\mathbf x}}

\newcommand{\ff}{{\mathbf f}}
\newcommand{\Div}{{\rm div}}
\newcommand{\symnabla}[1]{{D(#1)}}
\newcommand{\stress}{{\boldsymbol\sigma}}
\newcommand{\bdf}{{\sigma}}
\newcommand{\eye}{{I}}
\newcommand{\vertiii}[1]{{\left\vert\kern-0.25ex\left\vert\kern-0.25ex\left\vert #1 
    \right\vert\kern-0.25ex\right\vert\kern-0.25ex\right\vert}}
\newcommand{\normal}{{\mathbf n}}
\newcommand{\curv}{{H}}

\title{A reduced 3D-0D fluid-structure interaction model of the aortic valve that includes leaflet curvature}
\author[$^1$]{Ivan Fumagalli}
\author[$^1$]{Luca Dede'}
\author[$^{1,2}$]{Alfio Quarteroni}
\affil[$^1$]{\small MOX, Dipartimento di Matematica, Politecnico di Milano, Italy
\protect\\\texttt {\{ivan.fumagalli,luca.dede,alfio.quarteroni\}@polimi.it}}
\affil[$^2$]{\small Institute of Mathematics, \'Ecole Polytechnique F\'ed\'erale de Lausanne, Switzerland
}
\date{}

\begin{document}

\maketitle


\begin{abstract}
	{We introduce an innovative lumped-parameter model of the aortic valve, designed to efficiently simulate the impact of valve dynamics on blood flow.}
	Our reduced model includes the elastic effects associated {with} the leaflets' curvature and the stress exchanged with the blood flow.
	The introduction of a lumped-parameter model based on momentum balance entails an easier calibration of the model parameters: phenomenological-based models, on the other hand, typically have numerous parameters.
	This model is coupled to 3D Navier-Stokes equations describing the blood flow, where the moving valve leaflets are immersed in the fluid domain by a resistive method.
	A stabilized finite element method with a BDF time scheme is adopted for the discretization of the coupled problem, and the computational results show the suitability of the system in representing the {leaflet} motion, the blood flow in the ascending aorta, and the pressure jump across the leaflets.
	Both physiological and stenotic configurations are investigated, and we analyze the effects of different treatments for the {leaflet} velocity on the blood flow.
\end{abstract}

\section{Introduction}

Cardiac valves allow maintaining unidirectional blood flow in the heart and circulatory system. Their function helps in reducing blood stagnation by improving the chamber washout, and in directioning the ejection jets and the coherent vortex structures of the flow.
Because of the relevance of such components, several cardiac pathologies are directly related (or at least entail) valvular abnormal conditions, such as calcification, stenosis, regurgitation, and anatomical defects {of} the leaflets or the subvalvular apparatus: see, e.g., \cite{schoen2005cardiac,otto2008calcific,xanthos2011anatomic,el2018mitral}.

Due to the complexity of cardiac {valve} structure and their strong interplay with blood flow, {valve} modeling in computational hemodynamics has been developed with different levels of detail.
Several works consider a prescribed kinematics of the valves, introducing interface conditions that are expressed by analytical laws, as in \cite{doi:10.1098/rsfs.2010.0036,mittalMitralImposedKinematics,tagliabue2017fluid,tagliabue2017complex,ZINGARO2021380,zingaro2022geometric}, or derived from clinical measurements, as in \cite{vergaraBAV,chnafaLeftHeart,bonomiBAV,lvTurbulence,this2020pipeline}.
On the other hand, detailed mechanical models for the leaflets have also been proposed in the literature, possibly including
the inhomogeneities and the fibers in the leaflets (see, e.g., \cite{maromFibersForFluid,marom,kaiser2021design}) or a mechanical coupling with the subvalvular apparatus and the proximal vessels, such as in~\cite{chordae,chordaeEffect,pulmonarytrunk,valsalva1,valsalva2}.
In order to couple such complex models with hemodynamics, the solution of a three-dimensional Fluid-Structure Interaction (FSI) problem is required.
A wide range of numerical methods have been employed to this aim (\cite{bookIheart,fumagalli2024novel}), either in a boundary-fitting setting or from a  Eulerian standpoint: the Arbitrary Lagrangian Eulerian scheme (\cite{chengAleFsiAv,zhangAleValve,espinoAleValve,xale,vergaraFSIAV}), the CUTFEM and XFEM methods (\cite{xfemFsiZerothick,xfemSurfaceEquation,xfemFsiThick2D_and_dg4Zonca,xfemFsiThick3D1,xfemFsiThick3D2,xfemFsiThick3D3,xfemFsiThin2D,xfemZonca, xfemZoncaDg}), the immersed boundary (\cite{origIb,wangIbFem,borazjaniCurvib,borazjaniCurvibFsiAv,geCurvibAvCalcification,griffithPeskinIb,griffithIb,votta,bazilevsHughesIb,hsuImmersogeom,zhangIbFem,nestola2019immersed}) and the fictitious domain approach (\cite{glowinskiFd,fdFsiRigidcontact,fdSurfLambdaGerbeau,gerbeauFd,bazilevsFd,bazilevsHughesFd,dehartFdAV,dehartFdExperim,morsiAleTrileaflet}), the chimera method (\cite{chimeraValves,chimeraExperimental,chimeraLvRigidValve}), and space-time finite elements (\cite{hughesSpacetime,tezduyarSpacetime,tezduyarSpacetimeValve}), to mention a few.
The common ground of all these methods is that they require a full 3D (or at least 2D) representation of the valve geometry and of its mechanics solver, thus entailing a significantly increased computational cost with respect to imposed-displacement hemodynamics.

{Conversely, in many clinical applications, the primary focus is on blood flow dynamics rather than the stresses and mechanical response of the valves themselves.
In this regard, aiming} 
at modeling the {valve} dynamics with little computational burden, while retaining its interaction with the blood flow, lumped-parameter models have been introduced.
Most of such models, e.g.~those proposed in \cite{korshi,blanco20103d,circulation0D}, typically account for the valve hemodynamics effects by means of a phenomenological relationship between the pressure jump across the leaflets and the flowrate passing through them.
However, since the parameters appearing in the equations seldom have a precisely quantifiable physical meaning, the calibration of the model may be quite cumbersome and highly dependent on the specific application of interest.
Other works derive their reduced model from a momentum balance at the leaflets:
up to the authors' knowledge, this approach has been first adopted by \cite{pedrizzetti0d}, where the inertia and stiffness of the leaflets are neglected, while in the more recent work by \cite{mittal0d} a linear ordinary differential equation is introduced for the valve opening coefficient.
However, the geometry of the leaflets plays a marginal role in the model, affecting only the valve's inertia.

In this paper, we introduce a novel lumped-parameter structure model for the aortic valve, with the aim of enriching the description of the valve dynamics {with respect to other 0D models in the literature} while preserving a low computational effort {compared to} fully 3D FSI systems.
We derive our {simplified} model from the balance of forces at the leaflet, relating the elasticity of the leaflets to their curvature.
{Although not entailing a locally accurate description of the leaflet mechanics,} this approach allows to {synthetically}
account for the specific valve geometry and to relate it directly with the {total force exerted by the flow on the leaflets.}

Based on this mechanical model, we set up a 3D-0D fluid-structure interaction (FSI) system modeling the interplay between the three-dimensional blood flow in the ascending aorta and the aortic valve dynamics.
Blood dynamics is described by incompressible Navier-Stokes equations, and the hemodynamics effect of the valve's kinematics are accounted for by the Resistive Immersed Implicit Surface method (RIIS) introduced by \cite{resistivo}.
This method is inspired by the Resistive Immersed Surface (RIS) method by \cite{risFernandez,risAstorino},
it is characterized by a negligible computational overhead cost in CFD simulations, and its suitability for the description of hemodynamics effects of
cardiac valves has been shown
both in physiological and pathological conditions: see \cite{menghini,zingaro2022geometric} and \cite{sam,hcm}, respectively.
{Moreover, the proposed reduced model has also been applied to assess the effects of pulmonary valve replacement on the hemodynamics of the proximal pulmonary arteries \cite{criseo2024computational}.}

This paper is organized as follows.
In \cref{sec:modmeth} we introduce the FSI mathematical model, made of the novel lumped-parameter structural model of the aortic valve, the blood flow equations including the RIIS representation of the leaflets, and the coupling between the two systems.
The numerical approximation of the reduced FSI problem and the scheme for its solution are described in \cref{sec:discr}.
Then, computational results are presented in \cref{sec:results}{, including} an analysis of the reconstruction of {leaflet} velocity, the investigation of physiological and pathological conditions, and the comparison with a well-known model available in the literature.
{These results demonstrate that the proposed model serves as a computationally efficient tool for capturing the effects of valve dynamics on blood flow dynamics and patterns, as well as hemodynamics indicators that are meaningful to address clinical questions.}

\section{Models and methods}\label{sec:modmeth}

We present a reduced model for the fluid-structure interaction between the blood flow in the aorta and the aortic valve leaflets.
In \cref{sec:fluid}, we introduce the fluid dynamics system, with the valve effects modeled by the resistive method of \cite{resistivo,sam}.
Then, a reduced structure model for valve dynamics is derived in \cref{sec:0D}, considering the external forces induced on the leaflets by the surrounding blood, and the FSI coupling is presented in \cref{sec:FSI}.

\subsection{Fluid model and RIIS method}\label{sec:fluid}

\begin{figure}
\centering
\includegraphics[width=0.5\textwidth]{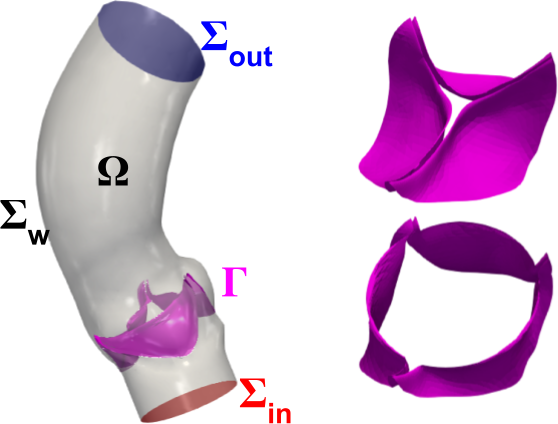}
\hfill
\begin{tikzpicture}
\draw[ultra thick] (0,0) -- node[right,rotate=90,anchor=south]{aortic wall} ++(0,5);
\draw[ultra thick, dashed] (5,0) -- node[right,rotate=90,anchor=north]{aorta centerline} ++(0,5);
\draw[thick] (0,0) to[out=0, in=-100] node[right]{$\Gamma_0$}(4.5,4);
\draw[thick] (0,0) to[out=0, in=-140] node[left]{$\Gamma_t$}(2,4.5);
\draw[->] (4,2.5) -- node[above]{$\dd_\Gamma(t,\xx)$} (1.5,3.5);
\draw[->] (2.5,0.7) -- node[right]{$\normal$} ++(-0.5,1);
\end{tikzpicture}
\caption{Computational domain and valve description. Left: the domain $\Omega$ with its boundaries and the immersed valve $\Gamma$ in purple; center: closed (above) and open (below) configuration of the aortic valve; right: schematic representation of a leaflet section and motion.}
\label{fig:domain}
\end{figure}

We model blood as incompressible and Newtonian,  with uniform density $\rho$ and viscosity $\mu$, and the domain $\Omega$ of interest is represented in \cref{fig:domain}.
The effects of the valve on the fluid dynamics are accounted for by the Resistive Immersed Implicit Surface (RIIS) method, introduced by \cite{resistivo} and employed by \cite{sam,hcm} in a clinical context.
This method, based on the Resistive Immersed Surface (RIS) approach proposed by \cite{risFernandez,gerbeauFd,risAstorino}, consists in the introduction of an additional penalty term in the fluid momentum equation, thus weakly imposing the kinematic condition at the surface representing the valve.

According to the RIIS method, the geometry of the moving valve $\Gamma_t$ is represented as a surface immersed in the fluid domain $\Omega$, implicitly described at each time $t$ by a level-set function $\varphi_t:\Omega\to\mathbb R$, as
\begin{equation}
\Gamma_t = \{\xx\in\Omega \colon \varphi_t(\xx)=\mathbf 0\}.
\end{equation}
The function $\varphi_t$ is assumed to be a signed distance function, namely to fulfill $|\nabla\varphi_t|=1$, for any $t$.
A smeared Dirac delta function $\delta_{t,\varepsilon}:\Omega\to[0,+\infty)$ is then introduced, to approximate the Dirac distribution -- rigorously, the codimension-1 Hausdorff measure --
with support on the surface $\Gamma_t$, as follows:
\begin{equation}
\delta_{t,\varepsilon}(\xx) =
\begin{cases}
\frac{1+\cos(\pi\varphi_t(\xx)/\varepsilon)}{2\varepsilon}
& \text{if }|\varphi_t(\xx)|\leq\varepsilon,\\
0
& \text{if }|\varphi_t(\xx)|>\varepsilon,
\end{cases}
\end{equation}
where the {half-amplitude} $\varepsilon$ is the smoothing parameter.

In these settings, the velocity $\uu$ and pressure $p$ of the blood satisfy the following formulation of the Navier-Stokes equations:
\begin{equation}\label{eq:NS}
\begin{cases}
\partial_t\uu + \rho\uu\cdot\nabla\uu - \nabla\cdot\stress + \frac{R}{\varepsilon}(\uu-\uu_\Gamma)\delta_{t,\varepsilon} = \mathbf 0
& \text{in }\Omega,\ t\in(0,T],\\
\nabla\cdot\uu = \mathbf 0
& \text{in }\Omega,\ t\in(0,T],\\
\uu = \mathbf 0
& \text{on }\Sigma_\text{w},\ t\in(0,T],\\
\stress\normal = p_\text{in}\normal,
& \text{on }\Sigma_\text{in},\ t\in(0,T],\\
\stress\normal = p_\text{out}\normal,
& \text{on }\Sigma_\text{out},\ t\in(0,T],\\
\uu = \mathbf 0
& \text{in }\Omega,\ t=0,
\end{cases}
\end{equation}
where $\stress=2\mu\symnabla{\uu}-p\eye = \mu(\nabla\uu+\nabla\uu^T)-p\eye$ is the fluid stress tensor, $R$ is the resistance of the RIIS term -- acting as a penalty parameter -- and $\uu_\Gamma$ is the velocity of the valve, which constitutes a data for the fluid problem and will be discussed in the following sections.
Regarding boundary conditions, $p_\text{in}, p_\text{out}$ are the pressure values imposed at the inflow and outflow boundaries $\Sigma_\text{in},\Sigma_\text{out}$, respectively, while the boundary $\Sigma_\text{w}$ represents the aortic wall.

\subsection{Lumped-parameters mechanical model}\label{sec:0D}

In order to provide the configuration and the velocity of the valve, represented by $\varphi_t$ and $\uu_\Gamma$ in the fluid problem \eqref{eq:NS}, a structural model would be required for the deformation of the surface $\Gamma_t$.
This section is devoted to the derivation of a reduced, lumped-parameters model realistically describing the main features of cardiac valve dynamics.
The approach differs from the one proposed by \cite{mittal0d} in that the elastic terms are related to the curvature of the leaflet, thus including additional geometrical information in the model.

Let $\dd_\Gamma:[0,T]\times\widehat\Gamma\to\mathbb R^3$ denote the displacement of the leaflet with respect to its reference configuration $\Gamma_0=\widehat\Gamma$, namely we can represent the current configuration $\Gamma_t$ as
\begin{equation}
\Gamma_t = \{\xx\in\mathbb R^3 \colon \xx=\dmapt(\widehat\xx)=\widehat{\xx}+\dd_\Gamma(t,\widehat{\xx}) \text{ for some }\widehat{\xx}\in\widehat\Gamma\},
\end{equation}
as schematically displayed in
\cref{fig:domain}.

We assume that at each time $t$, every point $\xx\in\Gamma_t$ of the leaflet is subject to an external force $\ff(t,\xx)$ due to the surrounding fluid and to an elastic force related to the leaflet curvature $\curv(\xx)$ -- both depending on the current configuration of $\Gamma_t$ described by $\dd_\Gamma(\xx)$ -- and that the valve motion can be affected by some damping effect.
Regarding the curvature-induced elastic force, we assume that it acts only normally to the surface, similarly to what happens in free-surface tension (see, e.g., \cite{surfaceTension,fumagalli2018free}).
Moreover, since it is generally observed that the resting state of the aortic valve is its closed configuration, we impose this elastic force to vanish on $\widehat\Gamma$.

According to these assumptions, a local force balance can be formulated as follows:
\begin{equation}\label{eq:balance}
\rho_\Gamma\ddot\xx + \beta\rho_\Gamma\dot\xx = \ff(t,\xx) - \gamma\left(\curv(\xx)-\widehat\curv(\widehat\xx)\right)\normal_\Gamma(\xx),
\end{equation}
where $\rho_\Gamma$
{is a parameter accounting for the inertia of the valve leaflets,} 
$\beta$
is a damping coefficient, $\gamma$
is an elasticity coefficient, and $\normal_\Gamma$ is the normal to the surface $\Gamma_t$.
The function $\widehat\curv(\widehat\xx)$ denotes the total curvature of the surface $\widehat\Gamma$ in the position $\widehat\xx=\dmaptinv(\xx)$ corresponding to $\xx$, that is the curvature of the resting configuration. 
{The parameter $\rho_\Gamma$ can be considered an \emph{effective surface density} of the valve, since it accounts for the mass $m_\Gamma$ of the leaflets through the relation
\[
m_\Gamma= \int_\Omega \rho_\Gamma\delta_{t,\varepsilon}\,d\Omega = \rho_\Gamma|\Gamma_t|.
\]
Since we also know that $m_\Gamma=\rho_\text{valve}\,\ell\,|\Gamma_t|$, where $\rho_\text{valve}$ is the actual density of the leaflets and $\ell$ their average thickness, we have $\rho_\Gamma=\rho_\text{valve}\,\ell$.
In the following, we consider the common assumption that $\rho_\text{valve}=\rho$ and adopt $\ell=0.25$mm as the {leaflet} thickness \cite{morganti2015patient,thubrikar2018aortic}. 
}
It is worth to point out that the parameter $\rho_\Gamma$ {could also} be tuned to different values -- without loss of generality -- in order to account for the inertia of the leaflets, which may be unknown in patient-specific settings or possibly be affected by added-mass effects (cf., e.g., \cite{addedMass}).

Aiming at reducing equation \eqref{eq:balance} to a 0D model, we assume that $\dd_\Gamma$ can be decomposed as
\begin{equation}\label{eq:decomposition}
\dd_\Gamma(t,\widehat\xx) = c(t)\dbar(\widehat\xx),
\end{equation}
where the spatial dependence of the displacement, represented by the function $\dbar:\widehat\Gamma\to\mathbb R^3$, is known, whilst the opening coefficient $c:[0,T]\to\mathbb R$ has to be modeled.
In these settings, the local balance \eqref{eq:balance} can be re-written as
\begin{equation}
\left(\ddot c(t) + \beta\dot c(t)\right)\rho_\Gamma\,\dbar(\dmaptinv(\xx)) 
= \ff(t,\xx) - \gamma\left(\curv(\xx)-\widehat\curv(\dmaptinv(\xx))\right)\normal_\Gamma(\xx).
\end{equation}
Taking the component along $\normal_\Gamma(\xx)$ and integrating over $\Gamma_t$, we obtain the following ordinary differential equation for $c$:
\begin{equation}\label{eq:0dmodel}\begin{gathered}
\ddot c + \beta\dot c +\eta(c,\ff) = 0,
\qquad
\text{where}
\\
\eta(c(t),\ff(t)) 
= \frac{\gamma\int_{\Gamma_t}\left(\curv(\xx)-\widehat\curv(\dmaptinv(\xx))\right)
- \int_{\Gamma_t}\ff(t,\xx)\cdot\normal_\Gamma(\xx)\,d\xx}{\int_{\Gamma_t}\rho_\Gamma\dbar(\dmaptinv(\xx))\cdot\normal_\Gamma(\xx)
},
\end{gathered}\end{equation}
where the dependence of $\eta$ on $c$ is implicit in its dependence from the curvature $H$: indeed, $H=-\text{div}_\Gamma\normal_\Gamma$
and the normal vector $\normal_\Gamma$ can be computed in terms of the derivatives of the function $\dmapt(\widehat\xx)=\widehat\xx+c(t)\dbar(\widehat\xx)$; a more precise definition of $\normal_\Gamma$ and $\curv$ will be introduced in \cref{sec:FSI}.
Equation \eqref{eq:0dmodel} can be completed by proper initial conditions on $c$ and $\dot c$, depending on the application of interest.

\subsection{Coupling of the fluid and structure models}\label{sec:FSI}

\begin{figure}
\centering
\includegraphics[width=0.5\textwidth]{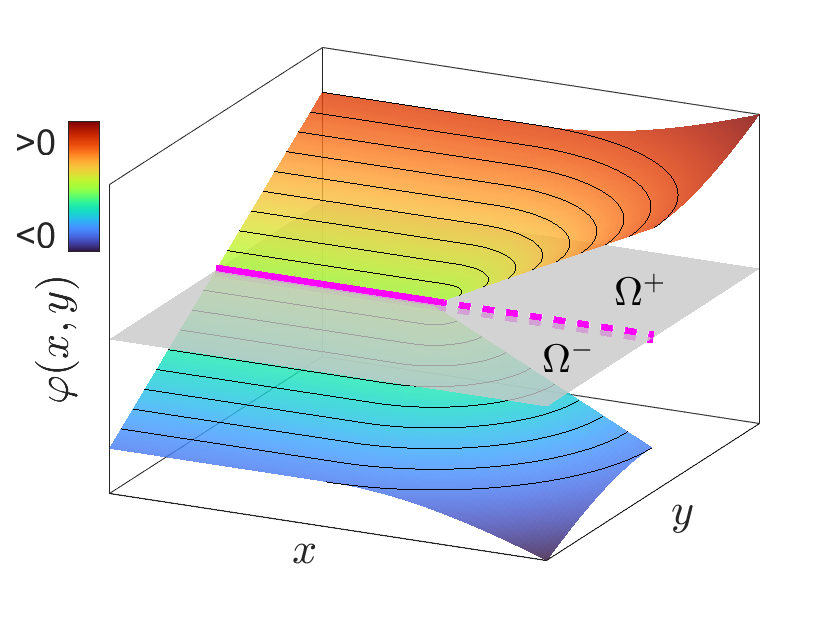}
\caption{Two-dimensional sketch of distance function $\varphi$ for a segment $\Gamma$ (solid magenta). In grey the plane $\varphi\equiv0$, with $\Omega^-$ and $\Omega^+$ separated by $\Gamma$ and the discontinuity line (dashed magenta).}
\label{fig:phidiscontinuous}
\end{figure}

We couple the 3D fluid model described in \cref{sec:fluid} and the 0D valve model introduced in \cref{sec:0D} to obtain a reduced FSI model:
the fluid-to-valve stress $\ff$ appearing in \eqref{eq:0dmodel} is computed from the former, while the latter provides the valve position and velocity.
To this aim, we introduce some additional notation related to the representation of the immersed surface $\Gamma_t$.
Being $\varphi_t$ a signed distance function, the domain $\Omega$ can be partitioned into two open sets
\begin{equation}
\Omega_t^+ = \{\xx\in\Omega \colon \varphi_t(\xx)>0\},
\quad
\Omega_t^- = \{\xx\in\Omega \colon \varphi_t(\xx)<0\}.
\end{equation}
Accordingly, any function $f$ defined over $\Omega$ can be decomposed as $f=f^++f^-$, where $f^\pm=f|_{\Omega^\pm}$.
\begin{rmrk}[Discontinuity of $\varphi_t$]
\label{rmrk:phidiscontinuous}
The definition of $\varphi_t$ that we employ, implemented in the Visualization Toolkit (VTK, \url{www.vtk.org}), yields that $\Gamma_t=\{\xx\in\Omega\colon\varphi_t=0\}$ is a \emph{subset} of the interface $\overline{\Omega_t^+}\cap\overline{\Omega_t^-}$ between $\Omega_t^-$ and $\Omega_t^+$.
Indeed, as schematically represented in \cref{fig:phidiscontinuous} for a 2D case with a segment $\Gamma$, such interface is partitioned into the actual leaflet $\Gamma$ and the line $(\overline{\Omega^+}\cap\overline{\Omega^-})\setminus\Gamma$ (a surface in 3D) where $\varphi$ jumps from negative to positive values.
\end{rmrk}

The function $\varphi_t$ allows to define $\widetilde{\normal}_\Gamma$ and $\widetilde{\curv}$, that are the extensions to the whole domain $\Omega$ of the surface normal $\normal_\Gamma$ and its curvature $\curv$, respectively:
\footnote{
In view of \cref{rmrk:phidiscontinuous}, the derivatives appearing in \eqref{eq:nhriis} are  computed in $\Omega^-$ and $\Omega^+$ separately, so that no contribution actually arises from the discontinuity of $\varphi_t$.
}
\begin{equation}\label{eq:nhriis}\begin{aligned}
\widetilde{\normal}_\Gamma &= \frac{\nabla\varphi_t}{|\nabla\varphi_t|},
\qquad
\widetilde{\curv} &= -\Div\,\widetilde{\normal}_\Gamma = -\frac{\Delta\varphi_t}{|\nabla\varphi_t|}+\frac{\nabla^2\varphi_t \colon (\nabla\varphi_t\otimes\nabla\varphi_t)}{|\nabla\varphi_t|^3},
\end{aligned}\end{equation}
with $\nabla^2\varphi_t$ denoting the Hessian matrix of $\varphi_t$.
The quantities $\widetilde{\normal}_\Gamma$ and $\widetilde{\curv}$ are actual extensions of the normal vector and curvature, since $\widetilde{\normal}_\Gamma|_\Gamma = \normal_\Gamma, \widetilde{\curv}|_\Gamma = \curv$ (cf., e.g., \cite{delfour2011shapes}).
We remark that $\widetilde{\normal}_\Gamma$ is such that it does not change its verse when passing through $\Gamma_t$.

\begin{rmrk}[Normalization]
In the definitions \eqref{eq:nhriis}, we did not make the standard assumption that $|\nabla\varphi_t|\equiv 1$.
Indeed, although such an assumption holds in the neighborhood of internal points of $\Gamma_t$, its validity is broken near $\partial\Gamma_t$, where $\varphi_t$ is not continuous.
Moreover, this definition of $\widetilde\normal_\Gamma$ ensures that the normal has unit magnitude also at the discrete level.
\end{rmrk}

Regarding the RIIS description of the surface, a definition of the surface velocity $\uu_\Gamma$ is required.
Based on the decomposition \eqref{eq:decomposition} of the displacement $\dd_\Gamma$, we provide the following definition:
\begin{equation}\label{eq:uGammaCont}
\uu_\Gamma(t,\xx) = \dot{c}(t)\widetilde\dbar(\xx),
\end{equation}
where $\widetilde\dbar:\Omega\to\mathbb R^3$ is the closest-point extension of $\dbar:\widehat\Gamma\to\mathbb R^3$.

The forces exerted by the fluid on the valve are related to the stress jump across $\Gamma_t$, thus
\begin{equation}
\ff = [\stress\normal_\Gamma]|_{\Gamma_t} = \stress^+|_{\Gamma_t}\normal_\Gamma - \stress^-|_{\Gamma_t}\normal_\Gamma.
\end{equation}
Considering the surface smearing introduced by the smooth Dirac delta $\delta_{\Gamma,\varepsilon}$ and the definitions \eqref{eq:nhriis}, the integral term related to $\ff$ that appears in \eqref{eq:0dmodel} can be approximated as follows:
\begin{align}\label{eq:fGamma}
\int_{\Gamma_t}\ff\cdot\normal_\Gamma &\simeq
\int_\Omega
\left(
\stress\widetilde\normal_\Gamma\cdot\widetilde\normal_\Gamma\,\delta^+_{\Gamma,\varepsilon}
- \stress\widetilde\normal_\Gamma\cdot\widetilde\normal_\Gamma\,\delta^-_{\Gamma,\varepsilon}
\right).
\end{align}

Analogously, the other integrals of \eqref{eq:0dmodel} can be approximated as follows:
\begin{equation}\label{eq:riisNormalCurvature}\begin{aligned}
\int_{\Gamma_t}\rho_\Gamma\left(\dbar\circ\dmaptinv\right)\cdot\normal_\Gamma &\simeq
\int_\Omega\rho_\Gamma\left(\dbar\circ\dmaptinv\right)\cdot\widetilde\normal_\Gamma\ \delta_{\Gamma,\varepsilon},
\\
-\gamma\int_{\Gamma_t}\left(\curv-\widehat\curv\circ\dmaptinv\right) &\simeq
-\gamma\int_\Omega
\left(
\widetilde\curv-\widehat{\widetilde\curv}\right)\delta_{\Gamma,\varepsilon},
\end{aligned}\end{equation}
with $\widehat{\widetilde\curv}$ denoting the RIIS representation of the pulled-back curvature $\widehat\curv\circ\dmaptinv$.

\begin{rmrk}[Transvalvular pressure jump]
Notice that, \linebreak since $|\widetilde\normal_\Gamma|\equiv 1$, if the strain component of the normal stress is assumed to be negligible {with respect to}~the pressure term, the integral force in \eqref{eq:fGamma} gets down to
\begin{equation}
\int_{\Gamma_t}\ff\cdot\normal_\Gamma \simeq
\int_\Omega
\left(
p\,\delta^+_{\Gamma,\varepsilon}
- p\,\delta^-_{\Gamma,\varepsilon}
\right),
\end{equation}
in accordance with other reduced models, such as those proposed by \cite{korshi,blanco20103d,mittal0d,pedrizzetti0d}, which are based on the pressure jump across the valve.
\end{rmrk}

\subsection{Numerical approximation}\label{sec:discr}

We present the space and time discretization of the coupled 3D-0D FSI model.
We introduce a uniform partition of the time interval $[0,T]$ with step-size $\Delta t$ and nodes $\{t^n=n\,\Delta t\}_{n=0}^N$.
Accordingly, the time-discrete counterparts of all quantities, evaluated at time $t^n$, will be denoted by a superscript $n$.
For the space discretization, we introduce a hexahedral mesh $\mathcal T_h$ for the domain $\Omega$, and the Finite Element (FE) space
\begin{equation}
X_h^r = \left\{v_h\in C^0(\overline\Omega) \colon v_h|_K \in \mathbb Q^r(K), \forall K\in\mathcal T_h\right\},
\end{equation}
where $\mathbb Q^r$ denotes the space of polynomials of degree $r$ {with respect to}~each space coordinate.
The velocity and pressure discrete spaces are thus defined as $V_h^r = \{\vv_h\in [X_h^r]^3\colon \vv_h=\mathbf 0 \text{ on }\Sigma_\text{w}\}$ and $Q_h^r=X_h^r$.

For the approximation of the fluid problem \eqref{eq:NS}, we adopt a semi-implicit BDF-FE scheme of order $\bdf$ as done by \cite{resistivo}, with the same polynomial degree $r$ for both $V_h^r$ and $Q_h^r$ and a SUPG-PSPG stabilization with VMS-inspired coefficients: cf.~\cite{fortivms,bazilevs2007vms}.

The resulting numerical method reads as follows:\\
Given $\uu_h^n\in V_h^r, n=0,\ldots,\bdf-1$, for each $n=\bdf,\ldots,N$, find $\uu_h^n\in V_h^r,p_h^n\in Q_h^r$ such that
\begin{align}\label{eq:fluidDiscr}
&\left(\rho\frac{\alpha_\bdf\uu_h^n-\uu_h^{n,\text{BDF}\bdf}}{\Delta t}, \vv_h\right)
+ \overline a^n(\uu_h^n,\vv_h)
+ c(\uu_h^{n,\bdf},\uu_h^n, \vv_h)
+ b(\vv_h,p_h^n) - b(\uu_h^n,q_h)
\\&\qquad
+\sum_{K\in\mathcal T_h}(\tau_\text{M}^{n,\bdf}\mathbf r_\text{M}^n(\uu_h^n,p_h^n), \rho\uu_h^{n,\bdf}\cdot\nabla\vv_h+\nabla q_h)_K
+\sum_{K\in\mathcal T_h}(\tau_\text{C}^{n,\bdf}r_\text{C}^n(\uu_h^n), \nabla\cdot\vv_h)_K = F(\vv_h)
\end{align}
for all $\vv_h\in V_h^r$ and $q_h\in Q_h^r$, where $(\cdot,\cdot)$ and $(\cdot,\cdot)_K$ denote the $L^2$ inner product over $\Omega$ and a mesh element $K$, respectively, and
\begin{align}
\overline a^n(\uu,\vv) &= \left(\mu\symnabla{\uu},\nabla\vv\right) + \left(\frac{R}{\varepsilon}\,\uu\,\delta_\varepsilon^n, \vv\right) ,
\\
b(\vv,q) &= -(\Div\vv,q) ,
\\
c(\mathbf w,\uu,\vv) &= \left(\mathbf w\cdot\nabla\uu,\vv\right) ,
\\
F(\vv) &= \int_{\Sigma_\text{in}} p_\text{in}\normal\cdot\vv + \int_{\Sigma_\text{out}}p_\text{out}\normal\cdot\vv 
- \left(\frac{R}{\varepsilon}\,\uu_{\Gamma,h}^n\,\delta_\varepsilon^n,\vv\right).
\end{align}
The BDF parameter $\alpha_\bdf$ and the velocities $\uu_h^{n,\text{BDF}\bdf}, \uu_h^{n,\bdf}$ depend on the order $\bdf$ of the BDF scheme (as in \cite{fortivms}), while
$\mathbf r_\text{M}^n, r_\text{C}^n, \tau_\text{M}^{n,\bdf}, \tau_\text{C}^{n,\bdf}$ are defined as
\begin{align}
\mathbf r_\text{M}^n(\uu_h^n,p_h^n) &= \rho\frac{\alpha_\bdf\uu_h^n-\uu_h^{n,\text{BDF}\bdf}}{\Delta t} - \mu\Delta\uu_h^n + \rho\uu_h^{n,\bdf}\cdot\nabla\uu_n^n
+ \nabla p_h^n +
\frac{R}{\varepsilon}\delta_\varepsilon^n(\uu_h^n-\uu_\Gamma^n),
\\
r_\text{C}^n(\uu_h^n) &= \nabla\cdot\uu_h^n,
\\
\tau_\text{M}^{n,\bdf} &= \left(\frac{\rho^2\alpha_\bdf^2}{\Delta t^2} + \rho^2\uu_h^{n,\bdf}\cdot \mathfrak{G}\uu_h^{n,\bdf} 
+ C_\text{r}\mu^2 \mathfrak{G}\colon \mathfrak{G} + \frac{R^2}{\varepsilon^2}(\delta_\varepsilon^n)^2\right)^{-1/2} ,
\\
\tau_\text{C}^{n,\bdf} &= \left(\tau_\text{M}^{n,\bdf}\mathfrak{g}\cdot\mathfrak{g} \right)^{-1}.
\end{align}
The quantities $\mathfrak{G}$ and $\mathbf{\mathfrak{g}}$ appearing above are the metric tensor and vector, depending on the element map $\mathbf M_K:\widehat{K}\to K$, for $K\in\mathcal T_h$, mapping the reference element $\widehat{K}$ to the current one $K$ (see, e.g., \cite{tezduyarSupg}).

\bigskip

Regarding the geometric quantities describing the valve, we hinge upon a FE description.
In particular, the discrete distance function is $\varphi^n_h\in X_h^{r'}$, with a polynomial degree $r'\geq2$ that is in general different from $r$.
Introducing the basis functions $\{\psi_\ell\}_{\ell=1}^{N_h^{r'}}$ spanning $X_h^{r'}$, the leaflet's extended normal and curvature are  defined as follows:
\begin{equation}\label{eq:basisFunNormalCurvature}\begin{aligned}
    \widetilde\normal_{\Gamma,h}^n &= \frac{\sum_{\ell=1}^{N_h^{r'}}\varphi_\ell^n\nabla\psi_\ell}{\left|\sum_{\ell=1}^{N_h^{r'}}\varphi_\ell^n\nabla\psi_\ell\right|},
    \\
    \widetilde H_{\Gamma,h}^n &= -\Div\,\widetilde\normal_{\Gamma,h}^n
    \\
    &= -\frac{\sum_{\ell=1}^{N_h^{r'}}\varphi_\ell^n\Delta\psi_\ell}{\left|\sum_{\ell=1}^{N_h^{r'}}\varphi_\ell^n\nabla\psi_\ell\right|}
    + \frac{\sum_{\ell,m,k=1}^{N_h^{r'}}\varphi_\ell^n\varphi_m^n\varphi_k^n\ \nabla^2\psi_\ell\ \colon \left(\nabla\psi_m\otimes\nabla\psi_k\right)}{\left|\sum_{\ell=1}^{N_h^{r'}}\varphi_\ell^n\nabla\psi_\ell\right|^3}
\end{aligned}\end{equation}
We point out that, since both these quantities appear in the valve model only as integrands of \eqref{eq:fGamma}-\eqref{eq:riisNormalCurvature}, we can use directly the expressions \eqref{eq:basisFunNormalCurvature}, without the need of a projection onto a finite element space.

Concerning the valve's kinematics, the discrete leaflet velocity is obtained from a first-order approximation of \eqref{eq:uGammaCont}:
\begin{equation}\label{eq:uGammaDiscr}
\uu_{\Gamma,h}^n = \frac{c^n - c^{n-1}}{\Delta t}\widetilde\dbar_h^n,
\end{equation}
while the solution of the ODE equation \eqref{eq:0dmodel} describing the valve dynamics is based on an explicit fourth-order Runge-Kutta method RK4 (cf.~\cite{bookForRK4}).

\begin{figure*}
\centering
\includegraphics[width=0.7\textwidth]{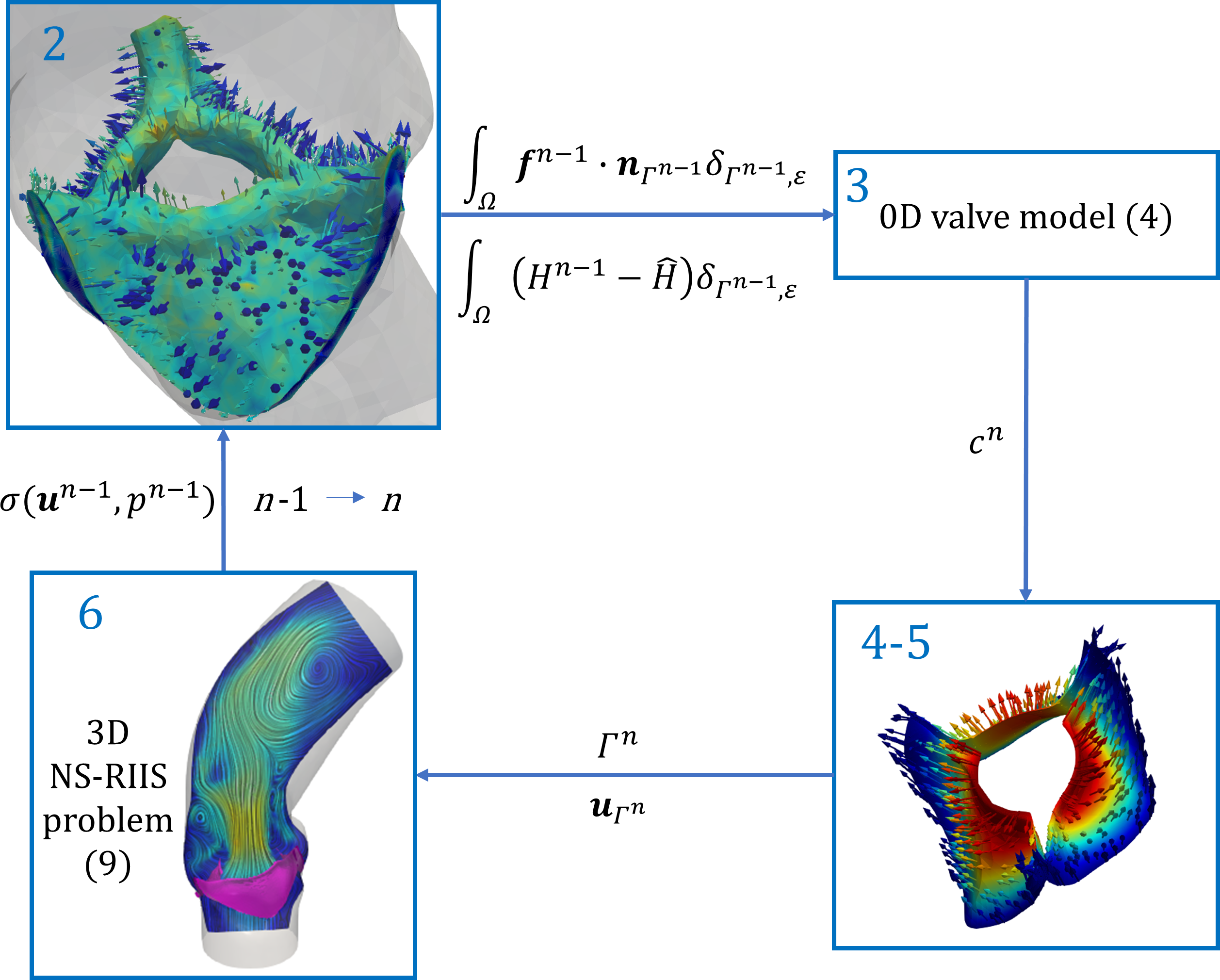}
\caption{Graphical representation of the staggered FSI solution scheme: the numbers of the panels correspond to the lines of Algorithm \ref{algo}. Panel 1: fluid-to-leaflet normal stress $\mathbf f\cdot\normal$ in the region $\{|\varphi_\Gamma|<\varepsilon\}$. Panel 3-4: {leaflet} velocity field $\mathbf u_\Gamma$. Panel 5: blood velocity $\uu$ on a slice.}
\label{fig:algo}
\end{figure*}

The fluid and structure models are weakly coupled at each time-step, as described in the following scheme, graphically displayed in \cref{fig:algo}:

\begin{algorithm}[H]
    \SetAlgoLined
    Given $\uu_h^n, p_h^n, c^n$ for $n=0,\dots,\bdf-1$, and computed the functions $\varphi^n,\widetilde\normal_\Gamma^n,\widetilde\curv^n$ corresponding to the surface $\Gamma^n$, for $n=0,\dots,\bdf-1$,
    \BlankLine
    \For{$n=\bdf$ \KwTo $N$}{
      Compute the integrals that make up \eqref{eq:0dmodel}, in terms of $\uu_h^{n-1},p_h^{n-1},\Gamma^{n-1},\varphi^{n-1}$\;
      Find $c^n$ by advancing the 0D equation \eqref{eq:0dmodel} with a step of RK4\;
      Move the immersed surface to its new configuration $\Gamma^n$ described by $\dd_\Gamma^n=c^n\dbar$ and compute $\uu_\Gamma^n = \frac{c^n-c^{n-1}}{\Delta t}\,\widetilde\dbar$\;
      Compute the new signed distance function $\varphi^n$ with respect to $\Gamma^n$ and
      assemble the normal and curvature fields $\widetilde\normal_\Gamma^n$ and $\widetilde\curv^n$;
      Find $(\uu_h^n,p_h^n)\in V_h^r\times Q_h^r$ by solving the linear problem \eqref{eq:fluidDiscr}.
    }
    \caption{Solution scheme for the 3D-0D FSI model}
    \label{algo}
    \end{algorithm}

This solution scheme has been implemented within life\textsuperscript{x} (cf.~{\cite{africa2022lifexcore,africa2024lifexcfd}, \url{https://lifex.gitlab.io/}}), a high-performance parallel C++ library for the solution of multi-physics problems based on the \texttt{deal.II} finite element core described by \cite{dealii}.

\section{Results and discussion}\label{sec:results}

We show the suitability of the proposed reduced 3D-0D FSI model in describing blood and valve dynamics in the ascending aorta.
Both the geometry of the domain $\Omega$ and of the closed valve leaflets $\widehat\Gamma$ are taken from Zygote (cf.~\cite{zygote2014}), an accurate model of the physiological heart derived from scan acquisitions.
{To define the open configuration} $\Gamma_\text{open}=\{\mathbf x = \widehat{\mathbf x} + \dbar(\widehat{\mathbf x}), \ \widehat{\mathbf x} \in\widehat\Gamma\}$ -- corresponding to an opening coefficient $c=1$ -- {we define $\dbar$ as proportional to the distance field $\widetilde{\dbar}$ connecting each point of a leaflet to the closest point to wall of the corresponding sinus of Valsalva.
Specifically, we progressively open the leaflets along $\widetilde{\dbar}$ until a physiological orifice area is attained.
The obtained valve configuration is shown in \cref{fig:domain}, bottom right, and has an orifice area of 2.78 cm\textsuperscript{2}, comparable with the values obtained in \cite{johnson2020thinner}.
}
{A possible drawback of this approach may be that the total area $|\Gamma_t|$ of the valve is not exactly constant throughout its motion, however the areas of the fully closed and fully open configuration differ by less than 1\%, and all intermediate configurations do not differ from them by more than 6\%.}

The domain is discretized by a hexahedral mesh of about 100K elements including artificial flow extensions at both inlet and outlet.
The elements size $h$ ranges from 2 mm in the flow extensions to 0.5 mm in the aortic root.
Blood velocity and pressure are both discretized with $\mathbb Q^1$ finite elements, and a BDF order $\bdf=1$, namely a semi-implicit Euler scheme, is chosen.
The other physical and numerical parameters of the system are reported in \cref{tab:param}.

\begin{figure}
\centering
\includegraphics[height=0.35\textwidth]{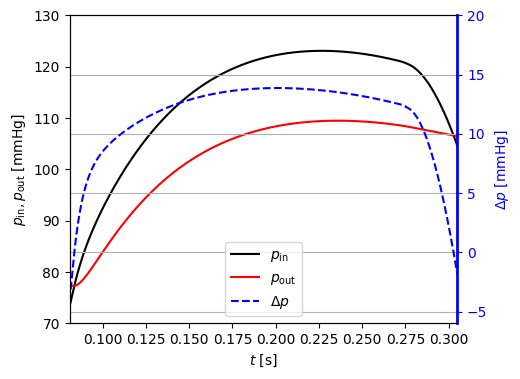}
\hfill
\includegraphics[height=0.35\textwidth]{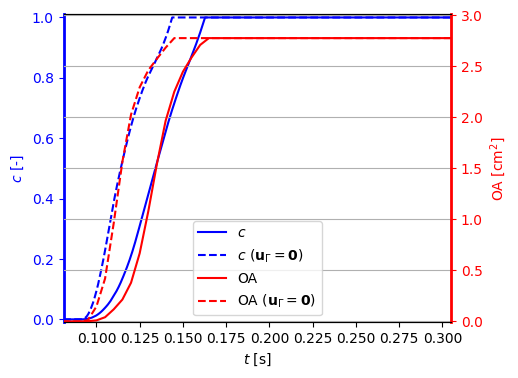}
\caption{Physiological pressure boundary conditions (left) and corresponding values of the valve's opening coefficient $c$ and orifice area (right).
The case $\uu_\Gamma=\mathbf 0$ is discussed in \cref{sec:quasistatic}.}
\label{fig:cSystole}
\end{figure}

Regarding boundary conditions at the inlet and outlet sections of the domain we impose the time-dependent normal stresses $p_\text{in}(t), p_\text{out}(t)$ displayed in \cref{fig:cSystole},
obtained from the lumped circulation presented by \cite{circulation0D} after proper calibration in order to be consistent with physiological pressures as reported in Wiggers diagrams (see, e.g., \cite{wiggersBook})
{and comparable with those employed in computational works such as \cite{bazilevsHughesFd,johnson2020thinner}.}

{
The choice of an effective calibration strategy is crucial to ensuring that the simulations accurately reflect blood flow physiology around the valve. Despite we did not conduct a systematic sensitivity analysis, we are aware of its importance in this context. Therefore, we explored a partial sensitivity analysis to assess the role of the following parameters in the valve's opening phase:
\begin{itemize}
    \item The damping parameter $\beta$ slows down the valve opening phase. Yet, to observe appreciable changes, the value of $\beta$ should be modified by at least one order of magnitude.
    \item Increasing $\gamma$ delays the opening phase and possibly prevents the valve from opening completely. A more detailed discussion is provided in \cref{sec:fullsystole}.
    \item An increase in the inertial parameter $\rho_\Gamma$ is associated with a slow down of the opening phase. A large increase of this parameter may further reduce the maximum attained orifice area. Further details are discused in \cref{sec:fullsystole}.
\end{itemize}
}

{All the simulations reported in the following were run in parallel on a 48-processor of CINECA's HPC cluster GALILEO100. On average, the wall time for the simulation of a full sysole was 5 hours. At each time step, the integration of the 3D stress terms on the leaflets and the solution of the reduced valve model required around one tenth of the computational time employed to solve the flow problem.
This shows how the proposed approach introduces little additional computational effort with respect to a purely fluid dynamics simulation.}

\begin{table}
\centering
\begin{tabular}{ccccccccc}
$\rho$ &
$\mu$ &
$R$ &
$\varepsilon$ &
{$\rho_\Gamma$} &
$\beta$ &
$\gamma$ &
$\Delta t$ \\
$\left[ \frac{\text{kg}}{\text{m}^3} \right]$ &
$\left[ \text{Pa}\ \text{s} \right]$ &
$\left[ \text{Pa}\ \text{s} \right]$ &
[m] &
{$\left[ \frac{\text{kg}}{\text{m}^2} \right]$} &
$\left[ \text{s}^{-1} \right]$ &
$\left[ \frac{\text{N}}{\text{m}} \right]$ &
[s]
\\ [2ex]%
\hline \\[-2ex]
        1060 &
$3.5\cdot10^{-3}$ &
$10^4$ &
$10^{-3}$ &
{$0.265$} &
$0.2$ &
{$3$} &
$2\cdot10^{-4}$
\end{tabular}
\caption{Physical and numerical parameters.}
\label{tab:param}
\end{table}

\subsection{Physiological valve opening}\label{sec:resPhysio}

\begin{figure*}
\centering
\begin{tabular}{ccccc}
a)
&
\includegraphics[height=0.3\textwidth]{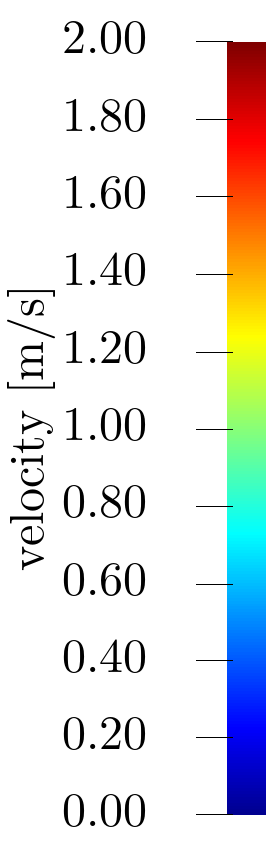}
&
\includegraphics[width=0.16\textwidth]{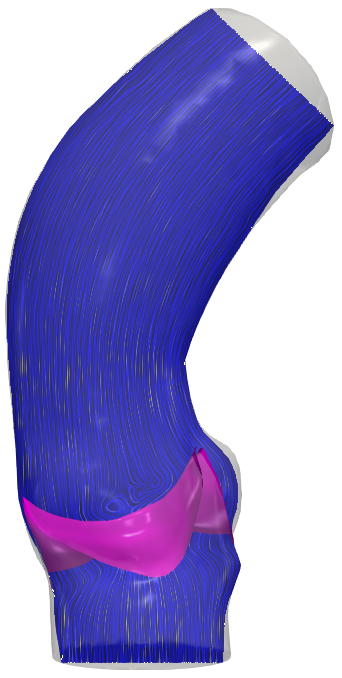}
&
\includegraphics[width=0.16\textwidth]{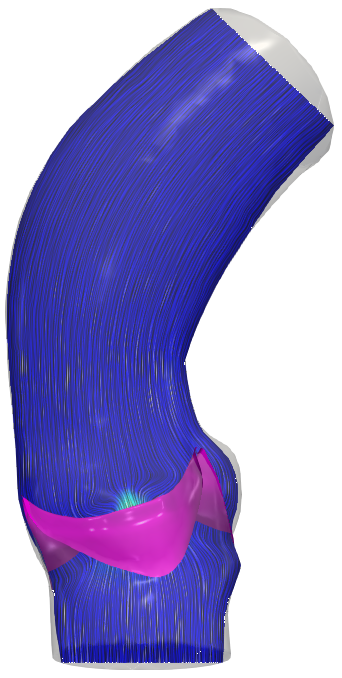}
&
\includegraphics[width=0.16\textwidth]{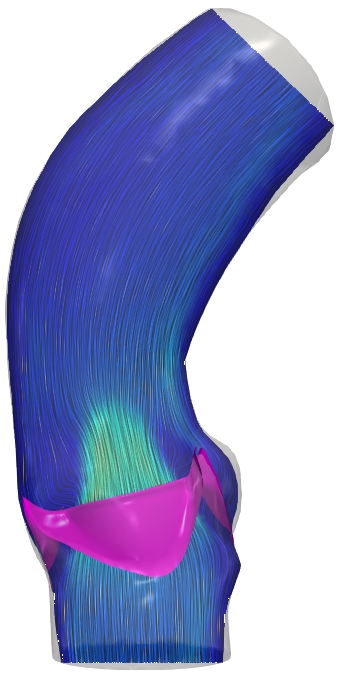}
\\
b)
&
\includegraphics[height=0.3\textwidth]{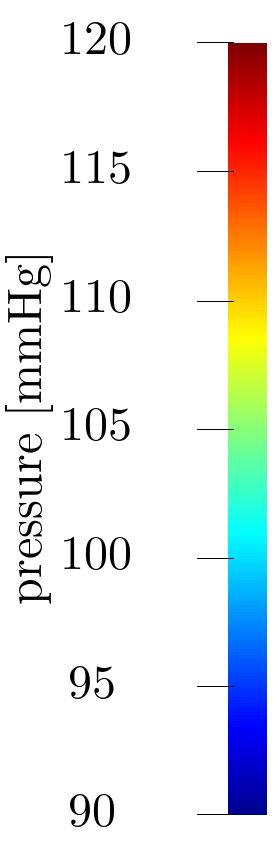}
&
\includegraphics[width=0.16\textwidth]{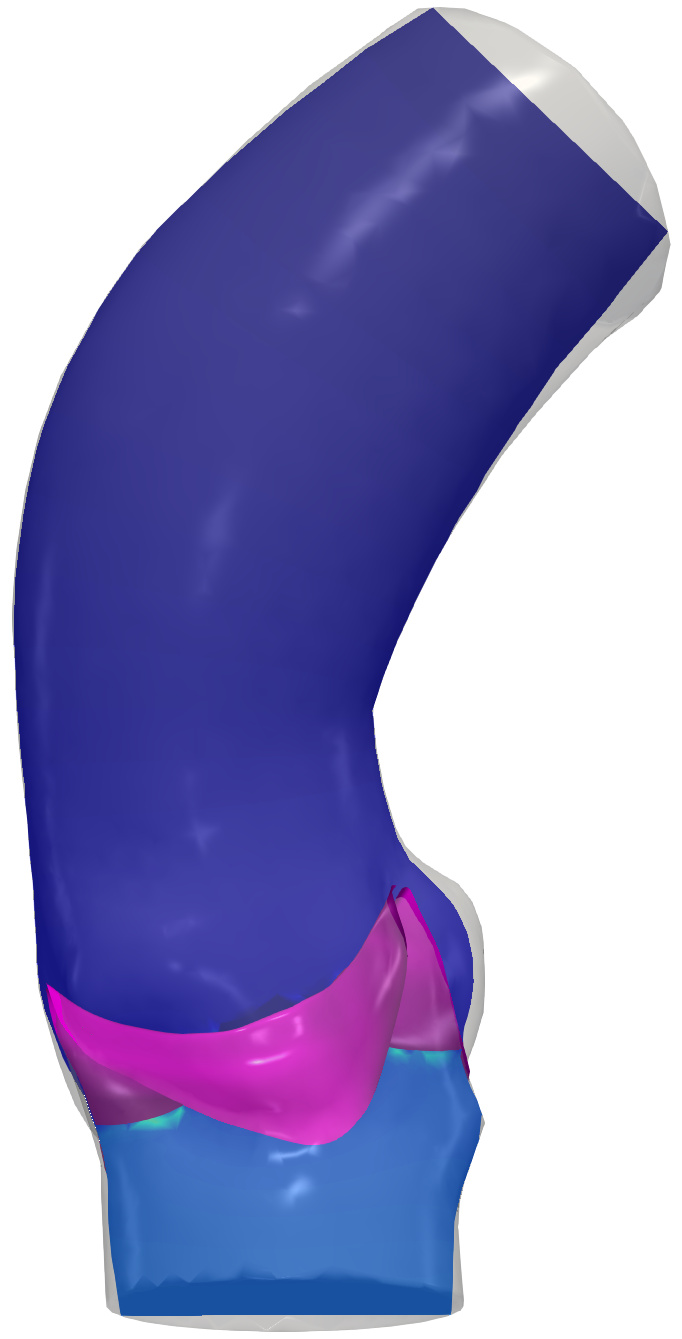}
&
\includegraphics[width=0.16\textwidth]{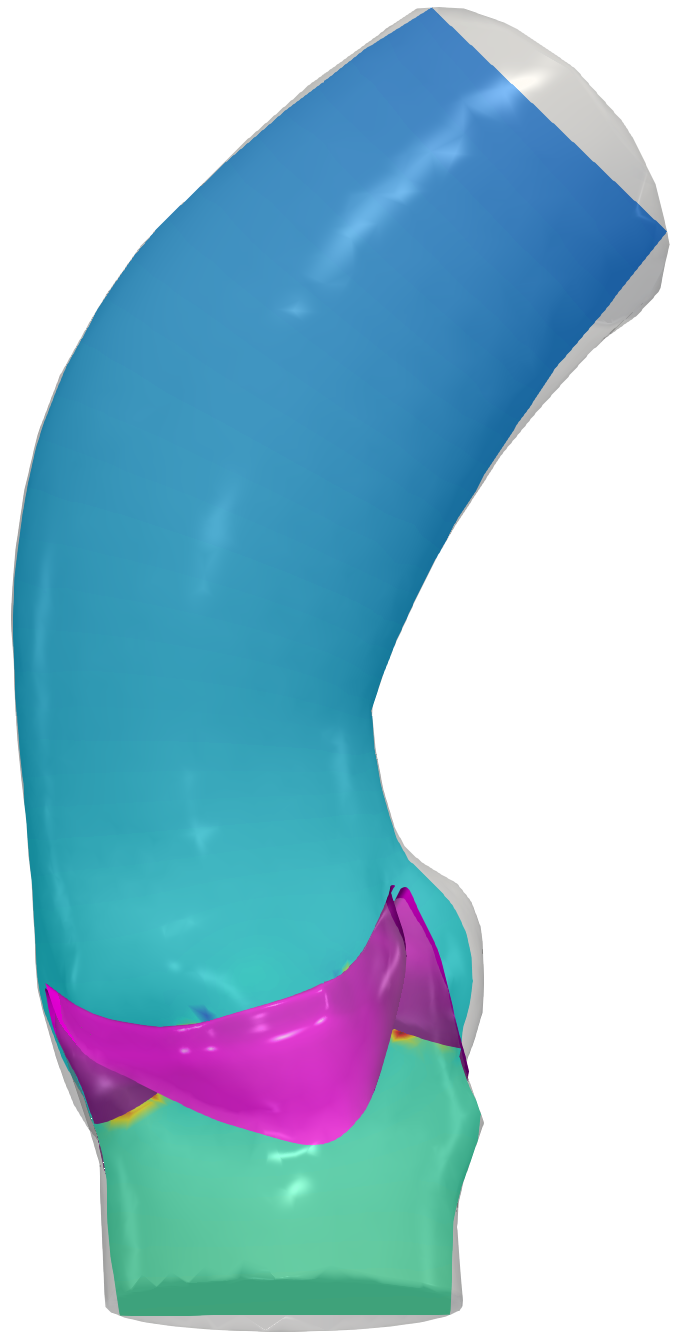}
&
\includegraphics[width=0.16\textwidth]{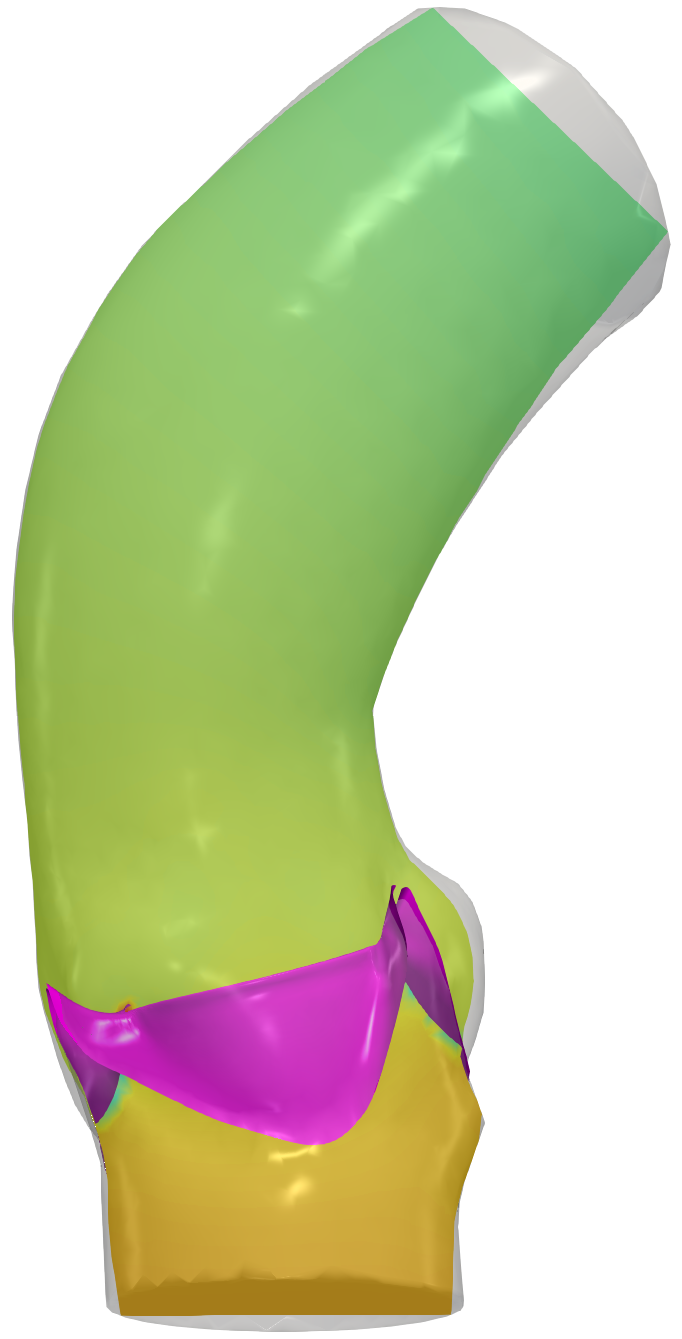}
\\
c)
& &
\includegraphics[width=0.23\textwidth]{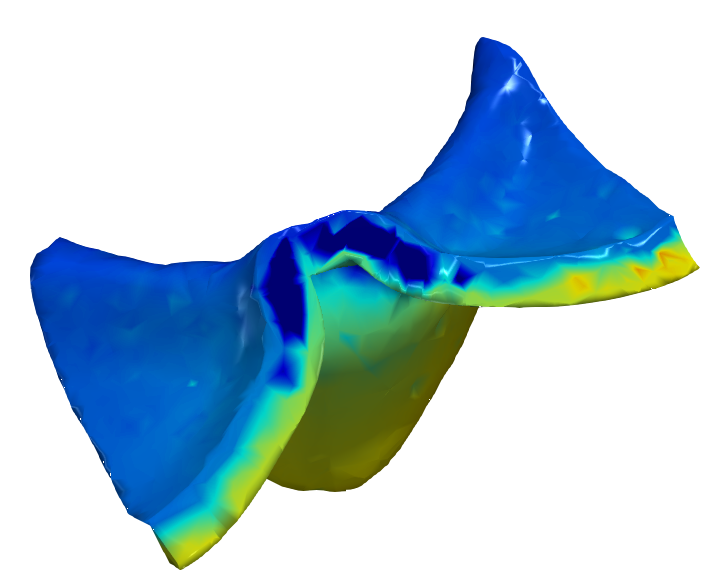}
&
\includegraphics[width=0.23\textwidth]{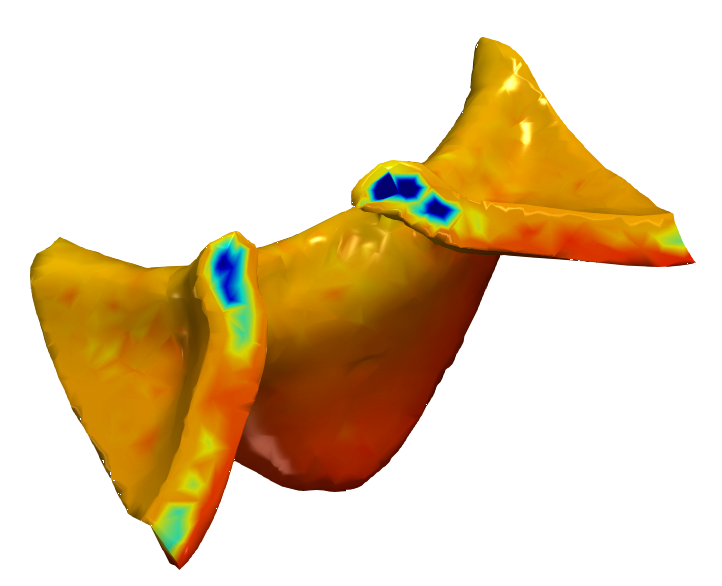}
&
\includegraphics[width=0.23\textwidth]{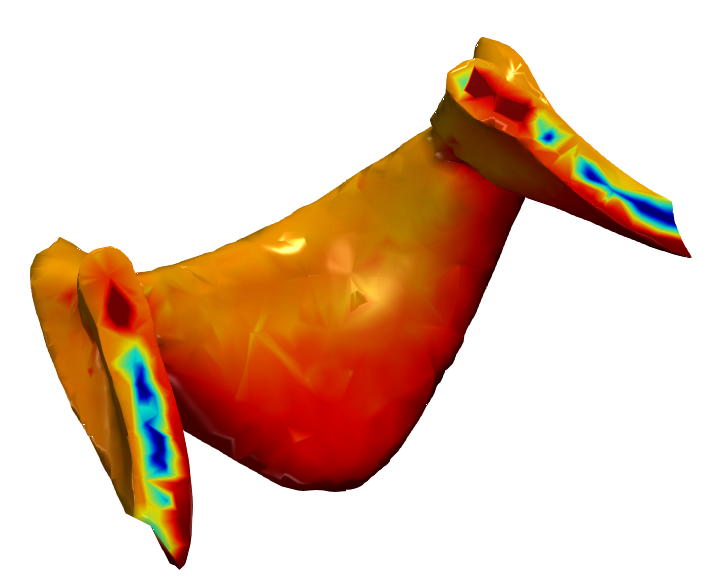}
\\[-1.4em]
& &
\includegraphics[width=0.1\textwidth]{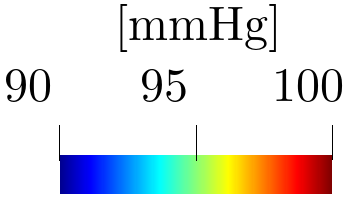}
&
\includegraphics[width=0.1\textwidth]{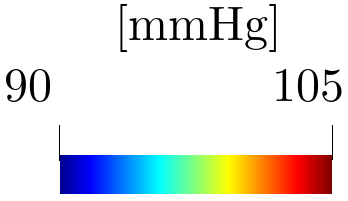}
&
\includegraphics[width=0.1\textwidth]{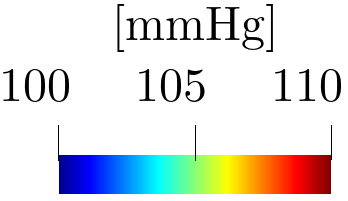}
\\
d)
&
\includegraphics[height=0.15\textwidth]{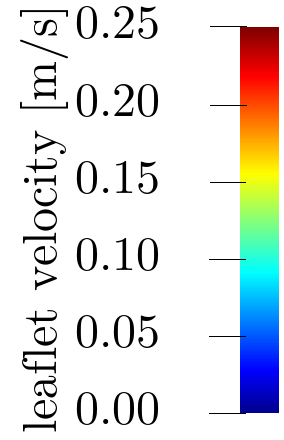}
&
\includegraphics[height=0.2\textwidth]{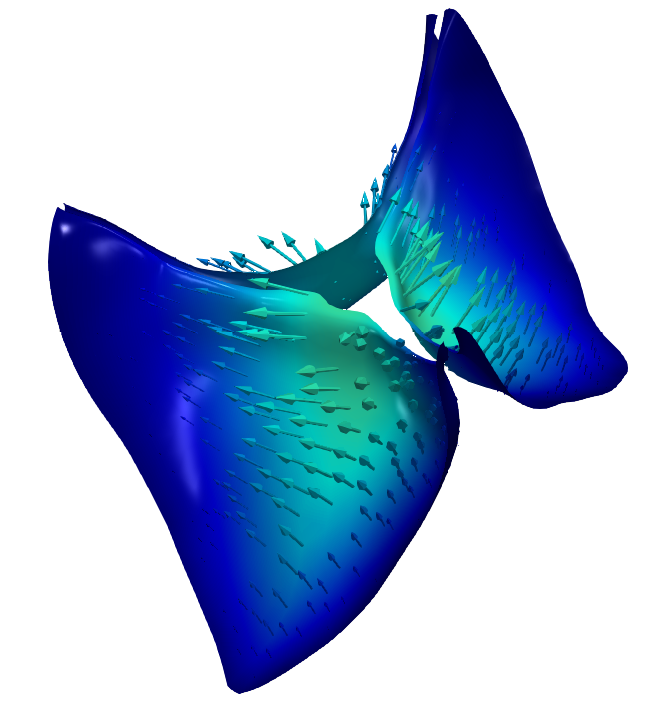}
&
\includegraphics[height=0.2\textwidth]{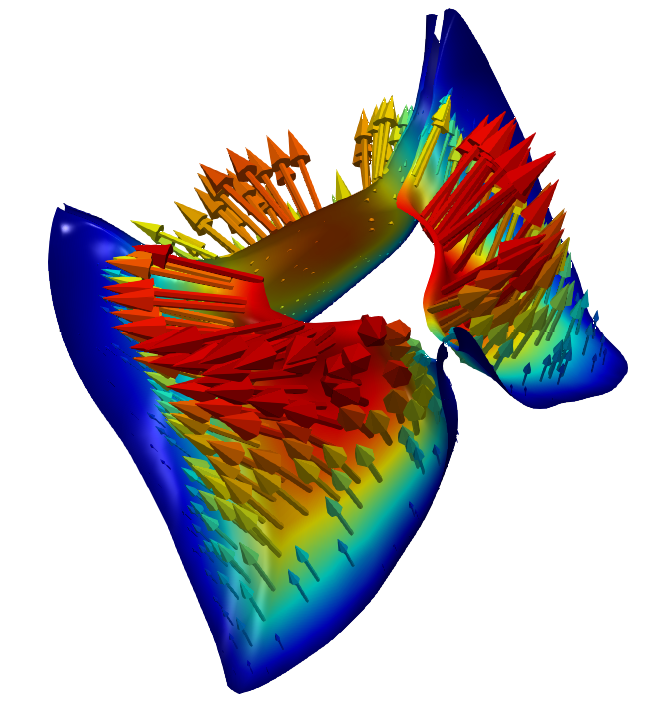}
&
\includegraphics[height=0.2\textwidth]{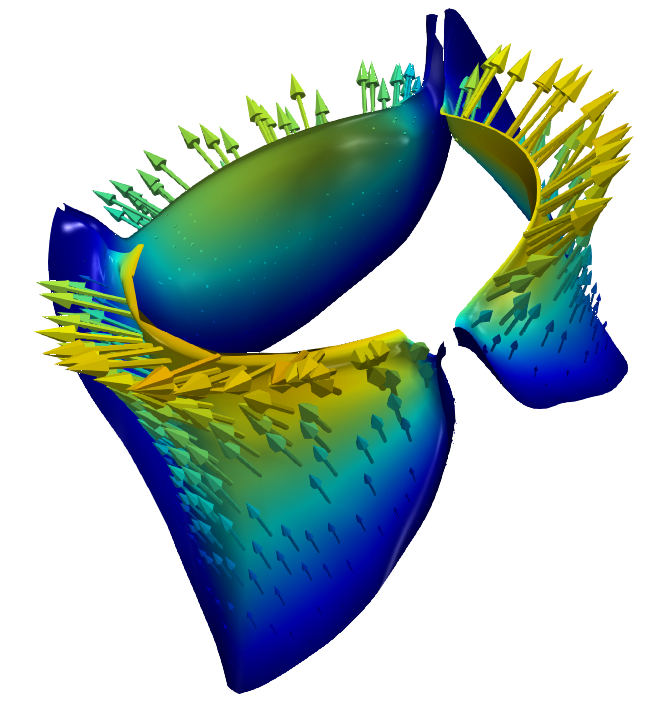}
\\
& & 
{$t=0.11$ s}
&
{$t=0.125$ s}
&
{$t=0.15$ s}
\end{tabular}
\caption{
Velocity (a) and pressure distribution in the domain (b) and within the {leaflet} region (c) under physiological pressure conditions. {Leaflet} velocity $\mathbf u_\Gamma$ in (d). The valve leaflets are colored in purple in (a)-(b).}
\label{fig:upSystole}
\end{figure*}

\begin{figure*}
\centering
\begin{tabular}{ccccc}
a)
&
\includegraphics[height=0.3\textwidth]{uLegend0to2.png}
&
\includegraphics[width=0.16\textwidth]{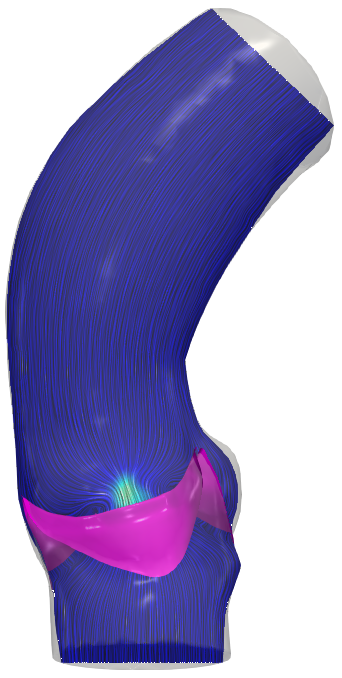}
&
\includegraphics[width=0.16\textwidth]{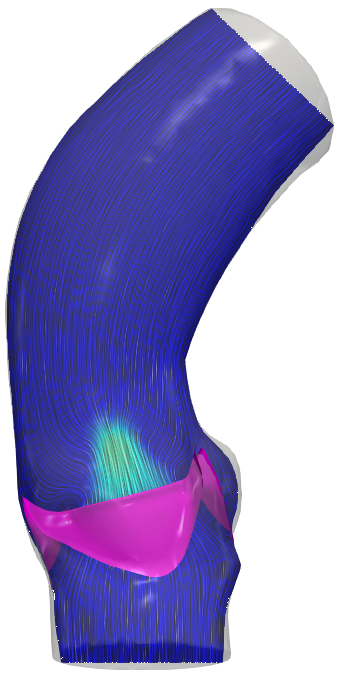}
&
\includegraphics[width=0.16\textwidth]{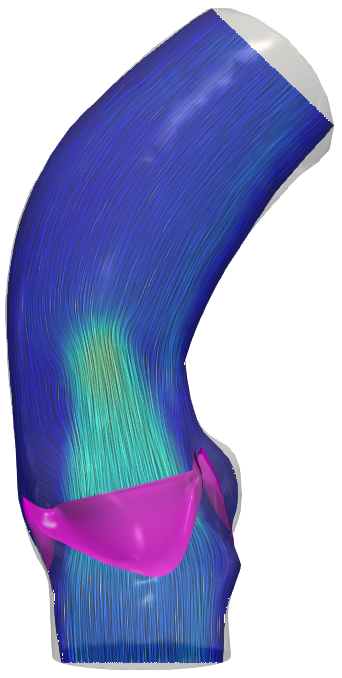}
\\
b)
&
\includegraphics[height=0.3\textwidth]{pLegend90to120.png}
&
\includegraphics[width=0.16\textwidth]{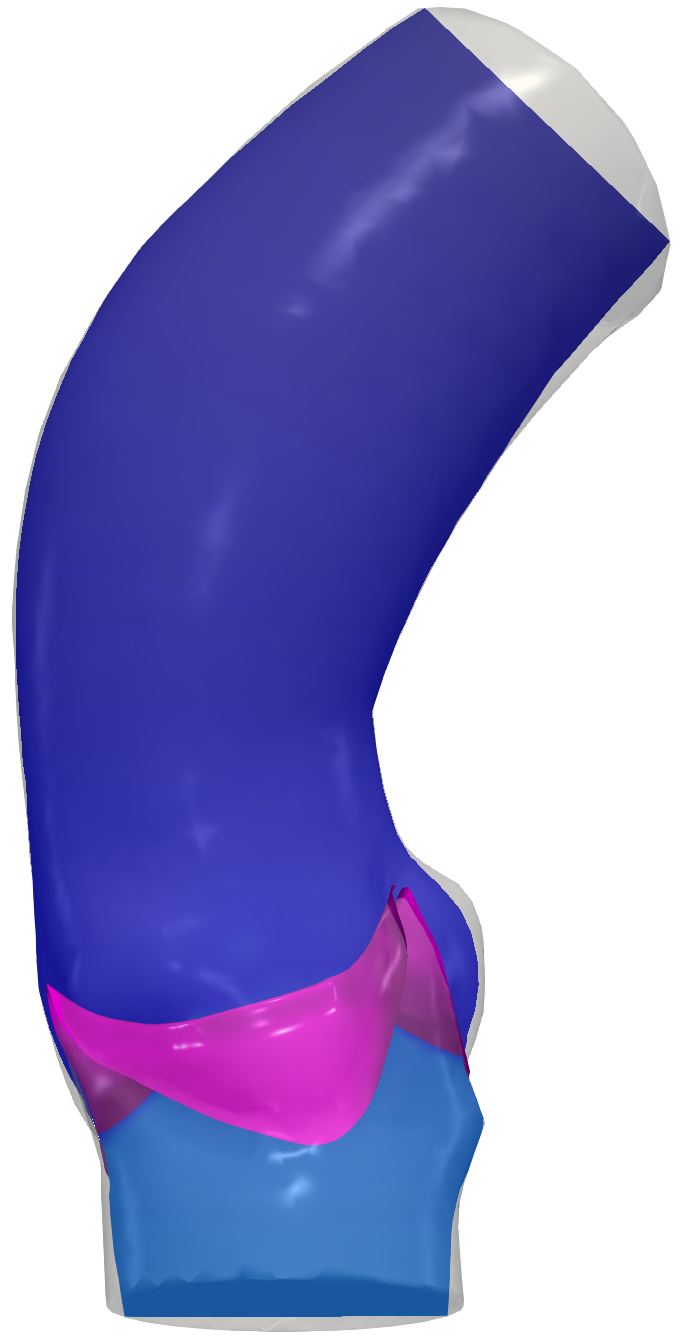}
&
\includegraphics[width=0.16\textwidth]{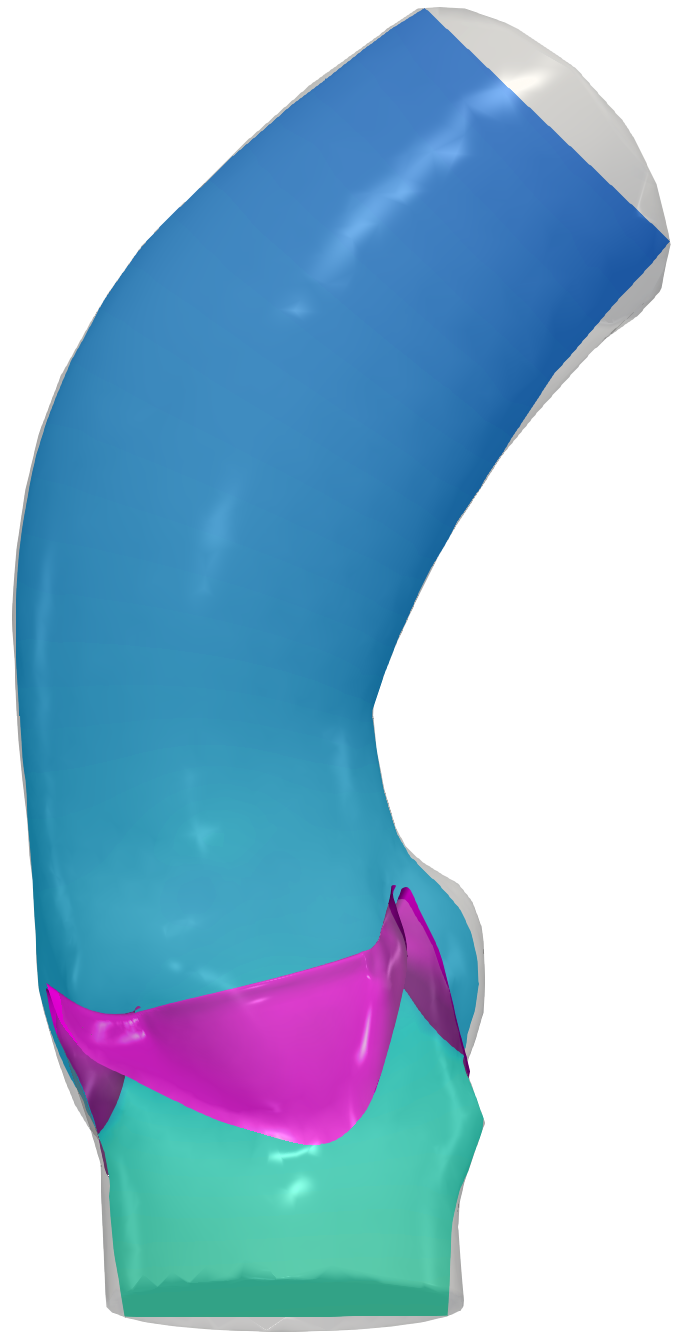}
&
\includegraphics[width=0.16\textwidth]{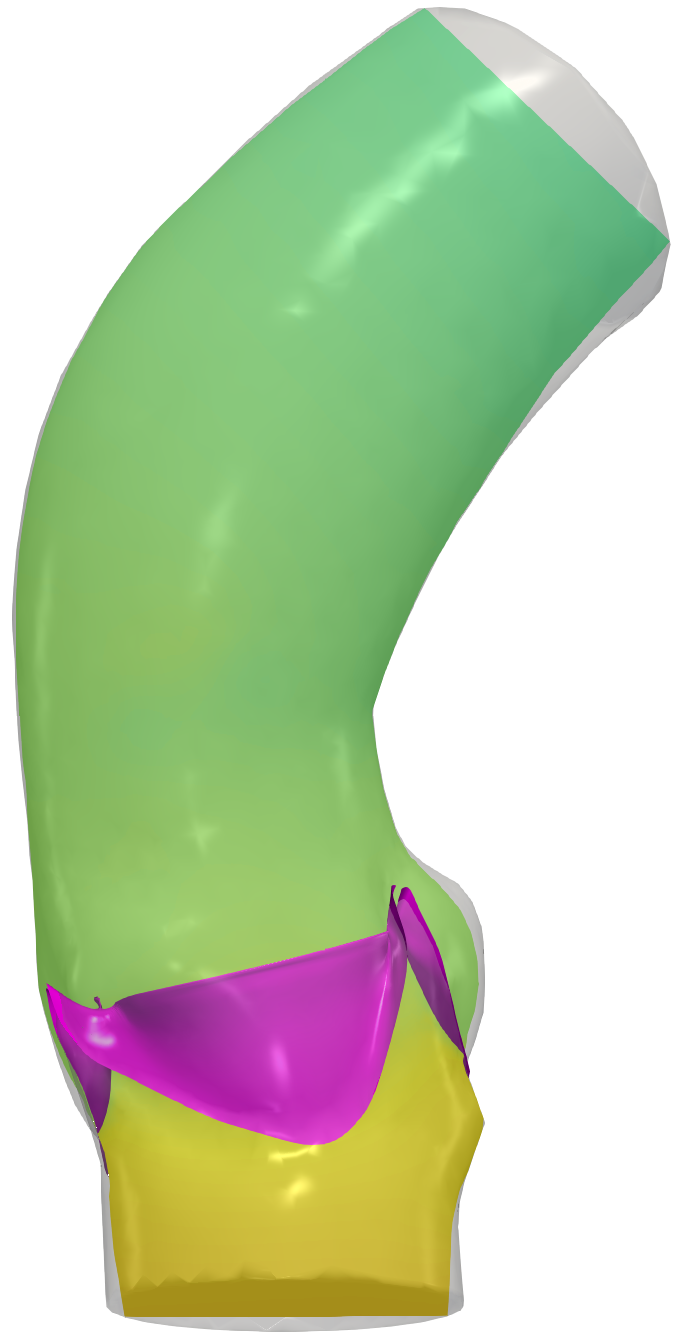}
\\
c)
& &
\includegraphics[width=0.23\textwidth]{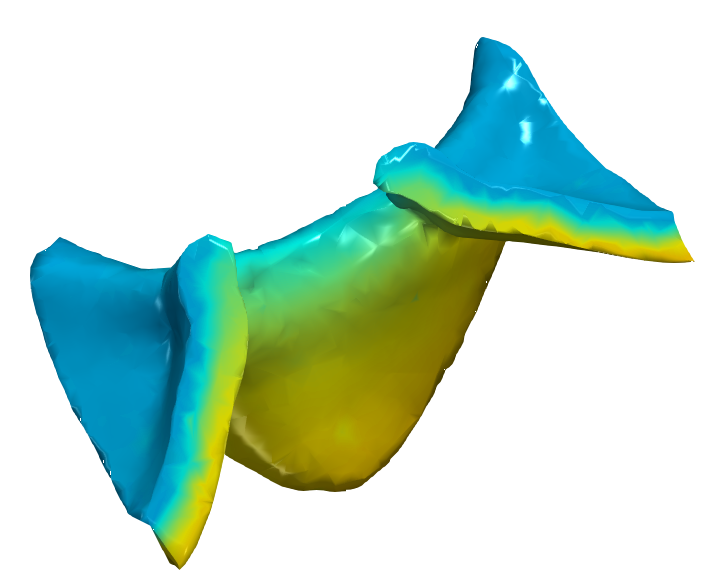}
&
\includegraphics[width=0.23\textwidth]{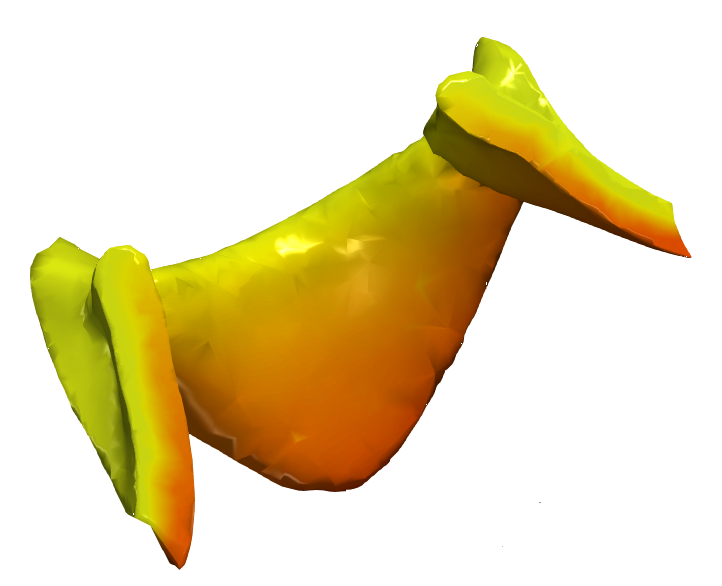}
&
\includegraphics[width=0.23\textwidth]{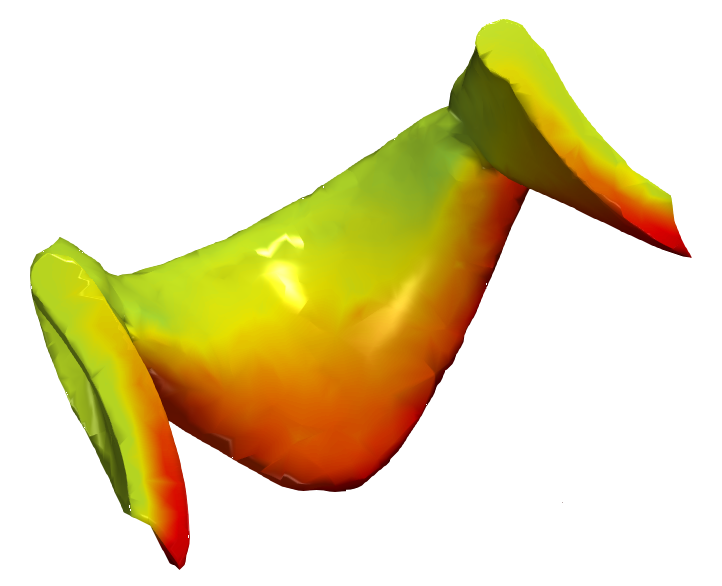}
\\[-1.4em]
& &
\includegraphics[width=0.1\textwidth]{pAVlegend90to100.png}
&
\includegraphics[width=0.1\textwidth]{pAVlegend90to105.png}
&
\includegraphics[width=0.1\textwidth]{pAVlegend100to110.png}
\\
& & 
{$t=0.11$ s}
&
{$t=0.125$ s}
&
{$t=0.15$ s}
\end{tabular}
\caption{
Quasi-static approach $\uu_\Gamma=\mathbf 0$. Velocity (a) and pressure distribution in the domain (b) and within the {leaflet} region (c) under physiological pressure conditions. The valve leaflets are colored in purple in (a)-(b).}
\label{fig:upSystoleNoUGamma}
\end{figure*}

We focus on the valve opening phase, {from $t=0.084$ s when the overall pressure difference $\Delta p=p_\text{in}-p_\text{out}$ between the inlet and the outlet is positive, up to its inversion occurring at $t=0.302$ s.}
As we can see from \cref{fig:cSystole}, the opening valve dynamics is characterized by different phases:
\begin{enumerate}
   \item[i)] The leaflets remain closed until a minimal transvalvular pressure jump of about 5 mmHg is developed.
   \item[ii)] Then, they rapidly open up to their fully open position, in a timespan of
   {72 ms,}
   in accordance with the measures of $76\pm 30$ ms reported by \cite{openingtime}.
   \item[iii)] In most part of the systole the valve remains in its fully open configuration, while the pressure jump progressively decreases.
\end{enumerate}

The evolution of the blood flow during this systolic ejection is reported in \cref{fig:upSystole}.
In the early stages of the simulation, while the valve is closed, the whole pressure gradient is concentrated across the valve.
Then, the opening of the valve is accompanied by a progressive development of the typical jet flow through the aortic orifice, and much smaller pressure differences can be observed.

In order to better examine the role of pressure in the valve dynamics, \cref{fig:upSystole}-c shows the pressure distribution in the $\varepsilon$-neighborhood of the leaflet, that is in the region where the RIIS term is active.
While the valve is closed, the whole pressure gradient develops within that region, showing the effectiveness of the RIIS method in providing an obstacle to the flow.
Then, while the valve opens, the pressure jump between the two sides of the leaflets is relatively small, but non-negligible gradients are present inside the RIIS region: this localized inhomogeneity allows to develop a nonzero leaflet velocity $\uu_\Gamma$ while preserving the incompressibility constraint of Navier-Stokes continuity equation.
Indeed, when the valve is in its fully open configuration, $\uu_\Gamma=\mathbf 0$ and pressure is essentially constant in the {whole $\varepsilon$-neigborhood of $\Gamma_t$.}

\subsection{Reconstruction of the {leaflet} velocity and quasi-static approach}\label{sec:quasistatic}

The leaflet velocity $\uu_\Gamma$ is provided by the reduced valve model \eqref{eq:uGammaDiscr}.
In this section, we assess its effect on the blood dynamics by comparing our results with those of the quasi-static approach adopted by \cite{resistivo}.
To this aim, a simulation in the same settings and boundary conditions of the previous section is run, the only difference being that $\uu_\Gamma=\mathbf 0$.
Resorting to \cref{fig:cSystole},
we can notice that the quasi-static approach {entails a faster opening phase (53 ms), with a larger opening velocity $\dot{c}$ especially at the beginning.} 
Moreover, comparing \cref{fig:upSystoleNoUGamma} with \cref{fig:upSystole}, a lower transvalvular pressure gradient can be observed at the early opening stages,
{as well as a faster developing jet in the aorta.}
These results can be motivated by observing that, in order to attain $\uu=\mathbf 0$ in the valve region, the continuous function $\uu$ must transition from the flow values to 0 in a surrounding boundary layer,
which thus artificially enlarges the effective obstacle that the leaflets represent to the flow{: as a consequence, the leaflets undergo a stronger push from the flow}.

We also compare our results with those of \cite{resistivo}, in terms of valve opening time.
It can be noticed that a much faster opening is observed in that reference (11 ms).
This difference is not only in the treatment of the surface velocity $\uu_\Gamma$, but also in the different valve model considered.
We can then state that the model presented in this work represents an improvement in terms of physiological representation of the aortic valve opening.
A more detailed comparison with such model is provided in \cref{sec:KS}.

\subsection{Full systole: physiological and stenotic valve}\label{sec:fullsystole}

\begin{figure}
\centering
\adjincludegraphics[width=0.14\textwidth, trim=0 0 18cm 0, clip=true]{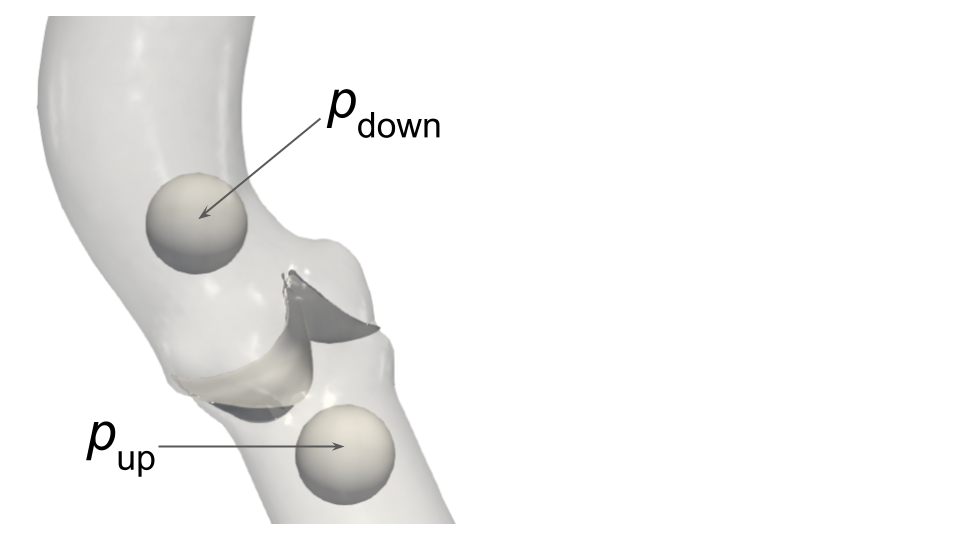}
\includegraphics[width=0.425\textwidth]
{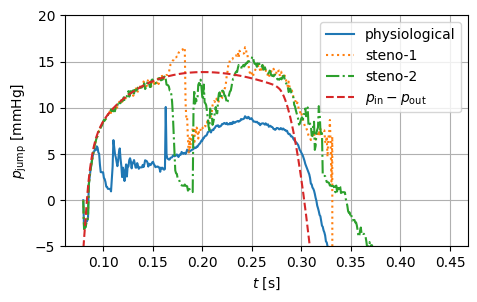}\hfill
\includegraphics[width=0.425\textwidth]
{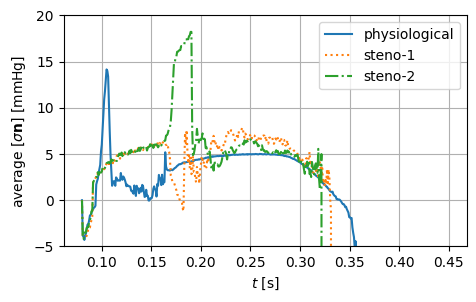}
\caption{Full systole: macroscopic pressure jump $p_\text{jump}=p_\text{up}-p_\text{down}$ (center) between two spherical control volumes (left) and transvalvular stress jump $\tfrac{1}{|\Gamma_t|}\int_{\Gamma_t}\mathbf f\cdot\normal_\Gamma$ (right). The overall pressure difference $\Delta p=p_\text{in}-p_\text{out}$ is reported, too, for comparison.}
\label{fig:pQFullSystole}
\end{figure}

We now employ the proposed reduced 3D-0D FSI model
to simulate a full systole, with the valve initially closed, namely $c(t=0)=0$.
In view of the discussion of \cref{sec:quasistatic}, we consider a non-zero leaflet velocity $\uu_\Gamma$, that is we do not adopt the quasi-static approach.
We are going to discuss our numerical results in a physiological case, and then we will introduce and investigate two different levels of aortic valve stenosis{, indicated as \emph{steno-1} and \emph{steno-2} in the following}.
{Specifically, case \emph{steno-1} corresponds to an increase of the elasticity coefficient $\gamma$ with respect to the physiological baseline -- modeling a stiffening of the valve -- and case \emph{steno-2} corresponds to an increase of the parameter $\rho_\Gamma$ with respect to \emph{steno-1} -- modeling an increase of the leaflets' inertia surrogating the added mass of calcifications.}
The values of $\gamma$ {and $\rho_\Gamma$} for the different cases are reported in \cref{tab:synth} together with the following synthetic indicators:
\begin{itemize}
    \item
$T_\text{open}$ is the time interval between the first time in which $c>0$ and the first one in which $c$ reaches its maximum value $c_\text{max}$ ($c_\text{max}=1$ in the physiological case);
analogously, $T_\text{close}$ is the time between the last local maximum of $c$ and the following instant in which $c=0$;
    \item
the aortic stenosis ratio $AS$ is based on the maximum orifice area $OA_\text{max}$:
$AS = 1-\frac{OA_\text{max}}{OA_\text{physio}}$ (cf.~\cite{mittal0d,OA}) with $OA_\text{physio}$ corresponding to the physiological case {$\gamma=3$ N/m;}    \item
$U_\text{peak}$ is the velocity attained by the aortic jet at the end of the opening phase;
    \item
$p_\text{jump, peak}$ is the macroscopic pressure jump $p_\text{jump}$ across the aortic root {at the time when $U_\text{peak}$ is attained}. This is computed as $p_\text{jump}=p_\text{up}-p_\text{down}$, where $p_\text{up},p_\text{down}$ are the average pressures in two small spheres upwind and downwind to the valve, respectively, as shown in \cref{fig:pQFullSystole}, left.
\end{itemize}

\begin{table}
\centering
\begin{tabular}{rl|ccc}
\multicolumn{2}{c|}{stenosis level} & physio. & {\emph{steno-1}} & {\emph{steno-2}}
\\
\hline
$\gamma$ & [N/m] & {3} & {15} & {15}
\\
{$\rho_\Gamma$} & {[kg/m\textsuperscript{2}]} & {0.265} & {0.265} & {0.276}
\\
$T_\text{open}$ & [ms] & {72} & {58} & {79}
\\
$T_\text{close}$ & [ms] & {32} & {35} & {21}
\\
$OA_\text{max}$ & [cm\textsuperscript{2}] & {2.78} & {1.79} & {1.06}
\\
$AS$ & [\%] & 0 & {36} & {62}
\\
$U_\text{peak}$ & [m/s] & {1.52}%
    & {2.26} & {2.17}
\\
$p_\text{jump, peak}$ & [mmHg] & {3.45}%
    & {8.41} & {12.52}
\end{tabular}
\caption{Synthetic indicators for valve stenosis. \label{tab:synth}}
\end{table}

\begin{figure}
\centering
\includegraphics[width=0.47\textwidth]{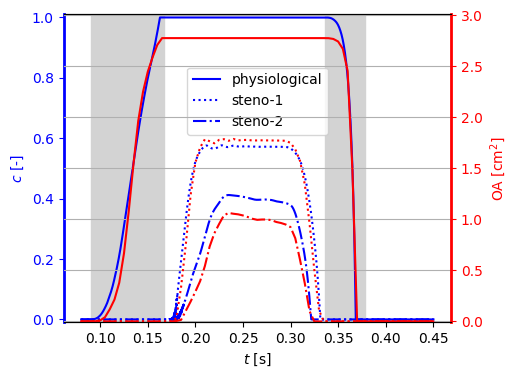}
\caption{Opening coefficient $c$ (left axis) and orifice area OA (right axis) under physiological pressure conditions, obtained with the curvature-based model in the case of a physiological{ valve and two degrees of aortic stenosis (\emph{steno-1}, \emph{steno-2}: see \cref{tab:synth}).}
The shaded areas correspond to average physiological opening and closing times as reported by \cite{openingtime}.}
\label{fig:cFullSystole}
\end{figure}

The baseline settings of the following discussion are those of the physiological case $\gamma={3}$ N/m.
As displayed in \cref{fig:cFullSystole}, after the opening phase discussed in \cref{sec:resPhysio}, the valve remains in its fully open position for {173}
ms, and then it closes in {32}
ms.
We point out that the duration of the closing phase lays within the physiological range of $42\pm 16$ ms reported by \cite{openingtime}, even though the calibration procedure considered only the \emph{opening} phase, thus supporting the aptness of the proposed valve model.

The time evolution of the pressure boundary conditions are reported in \cref{fig:pQFullSystole}, together with the cross-valve pressure jump $p_\text{jump}=p_\text{up}-p_\text{down}$.
We notice that $p_\text{jump}$ is definitely positive/negative during the opening/closing phase, whereas it remains below {10}
mmHg in the interval 
{$t\in[0.151,0.336]s$}
when the valve is fully open.
{These values are comparable, e.g., with the 5 and 15 mmHg of transvalvular pressure jump reported in \cite{bazilevsHughesFd} and \cite{johnson2020thinner}, respectively.
Moreover, the physiological value of the peak velocity $U_\text{peak}={1.52}$
m/s -- consistent with \cite{hsu2014fluid} --} 
confirms that the fully open state that we consider corresponds to a non-stenotic configuration.

We notice that the beginning of the closing phase at $t={0.368}$~s is delayed with respect to the inversion of the macroscopic pressure jump, occurring at $t={0.315}
$~s, and this delay is even larger than the closing time.
Such behavior, consistent with the valve modeling literature, is due to the inertia of both the blood flow and the valve, and it shows how the reconstruction of the {\em local} stress exchanged between the flow and the leaflets has a major impact on the valve dynamics.
Indeed, analyzing the time evolutions of \cref{fig:pQFullSystole}, we observe that,
in the interval $t\in{[0.151,0.336]}
$s where $c\equiv 1$, the stress jump $\tfrac{1}{|\Gamma_t|}\int_{\Gamma_t}\mathbf f\cdot\normal_\Gamma$ remains {between 3 and 5}
mmHg, keeping the valve open against the elastic forces; moreover, the change of sign in the stress term occurs at $t={0.34}
$~s, causing the abovementioned delay of the valve closing phase {with respect to}~the sign inversion of the pressure jump.
Furthermore, the average stress jump {remains significantly lower than $p_\text{jump}$ during almost all of the}
valve-opening phase: this can be seen as a confirmation of the common statement that cardiac valve leaflets (in physiological conditions) are basically \emph{transported} by the flow -- as done, e.g.~, in the purely kinematic model by \cite{collia2019simplified}.

\begin{figure*}
\centering
\begin{tabular}{m{1ex}ccccc}
a) &
\includegraphics[width=0.17\textwidth]{REV-uSLIC_kami78583_t0p15.png}
&
\includegraphics[width=0.17\textwidth]{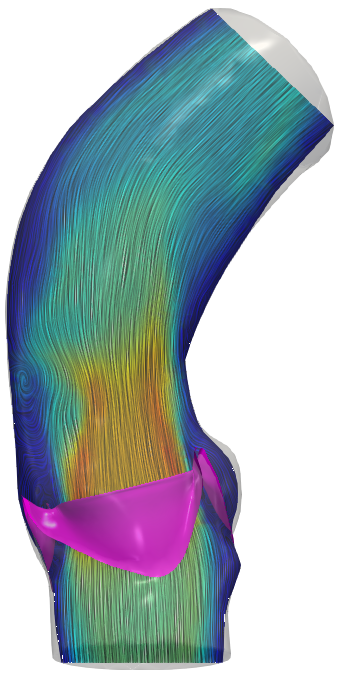}
&
\includegraphics[width=0.17\textwidth]{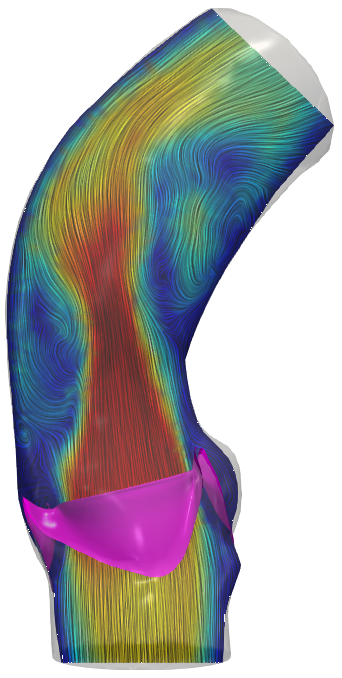}
&
\includegraphics[width=0.17\textwidth]{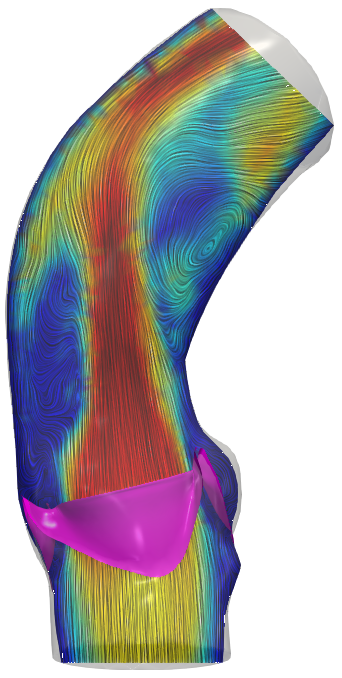}
&
\includegraphics[width=0.17\textwidth]{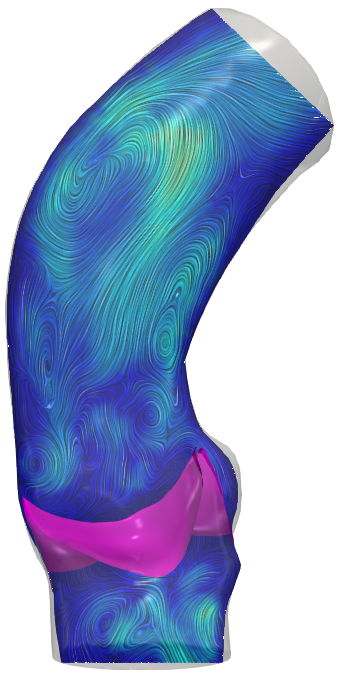}
\\
b) &
\includegraphics[width=0.17\textwidth]{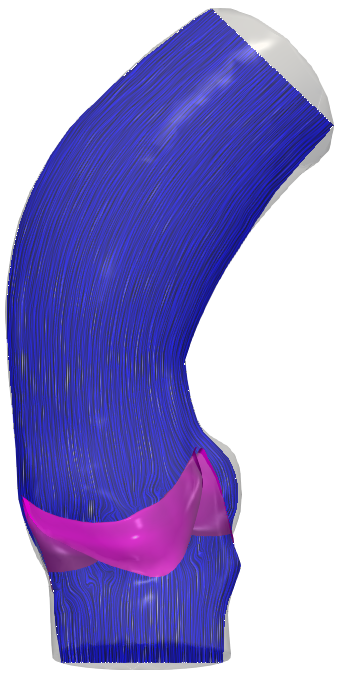}
&
\includegraphics[width=0.17\textwidth]{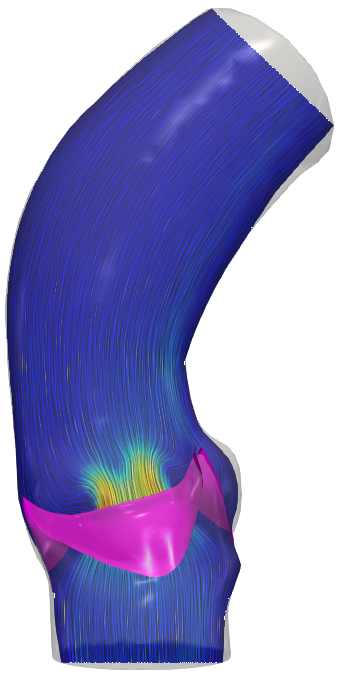}
&
\includegraphics[width=0.17\textwidth]{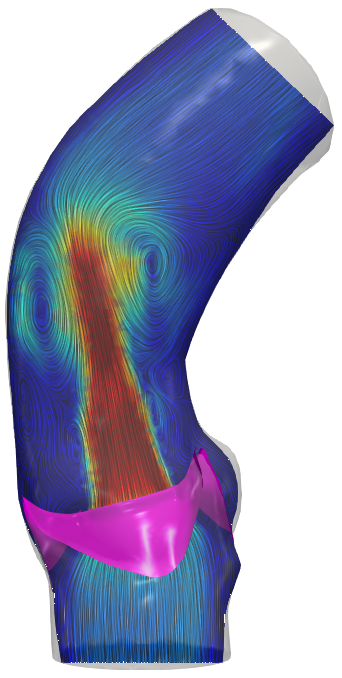}
&
\includegraphics[width=0.17\textwidth]{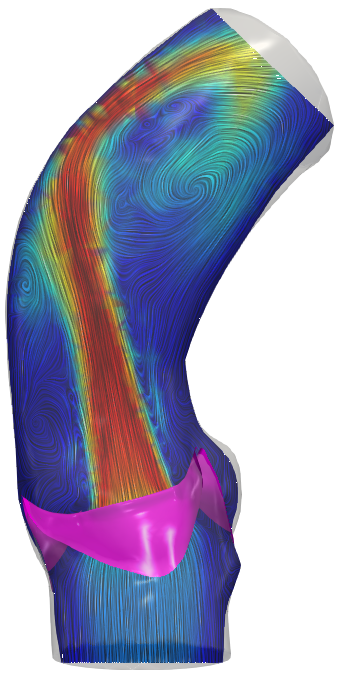}
&
\includegraphics[width=0.17\textwidth]{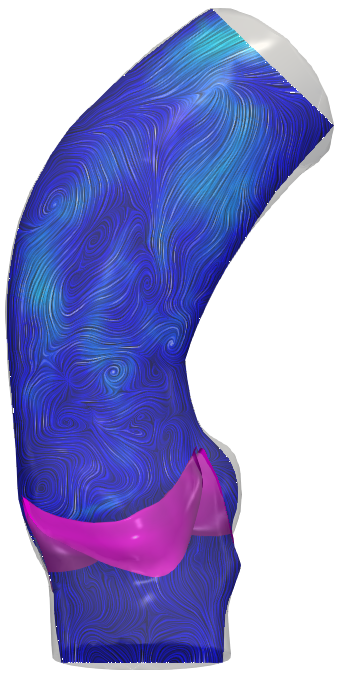}
\\
c) &
\includegraphics[width=0.17\textwidth]{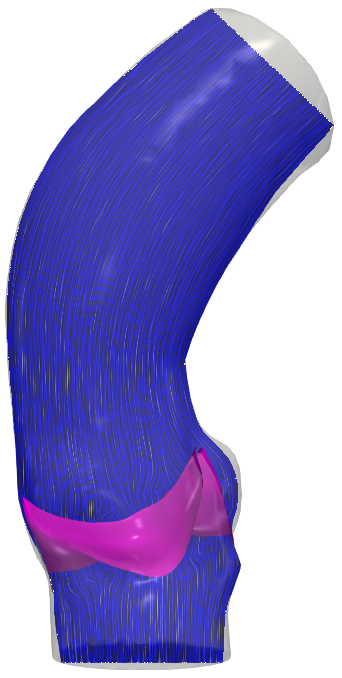}
&
\includegraphics[width=0.17\textwidth]{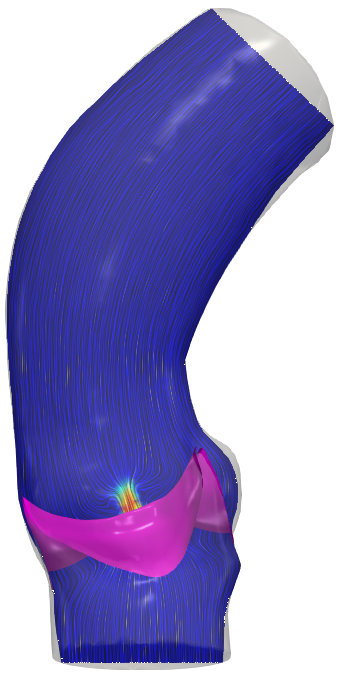}
&
\includegraphics[width=0.17\textwidth]{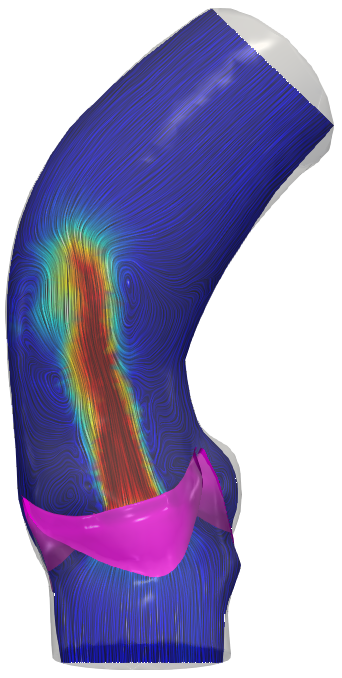}
&
\includegraphics[width=0.17\textwidth]{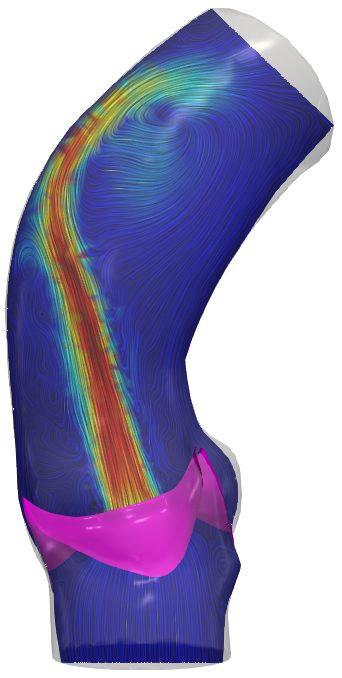}
&
\includegraphics[width=0.17\textwidth]{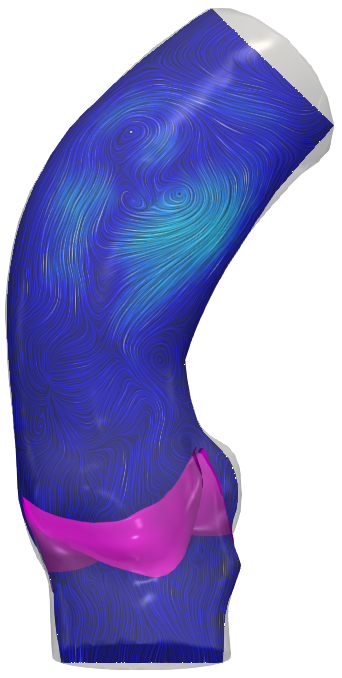}
\\
&{$t=0.15$ s} & {$t=0.20$ s} & {$t=0.25$ s} & {$t=0.30$ s} & {$t=0.40$ s} \\
\end{tabular}\\[1ex]
\ \qquad\qquad\includegraphics[width=0.4\textwidth]{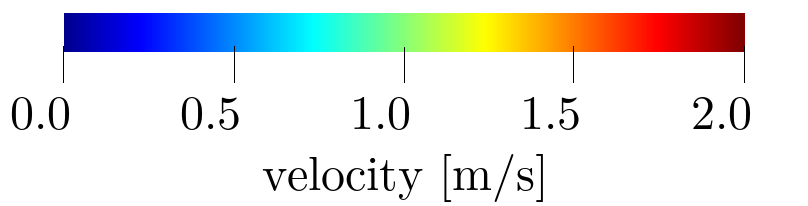}
\caption{
Velocity distribution on a longitudinal slice at different times:  a) physiological case, b) {case \emph{steno-1}}, c) {case \emph{steno-2} (see \cref{tab:synth})}.}
\label{fig:uSLIC}
\end{figure*}

\begin{figure*}
\centering
\begin{tabular}{m{1ex}ccccc}
a) &
\includegraphics[width=0.17\textwidth]{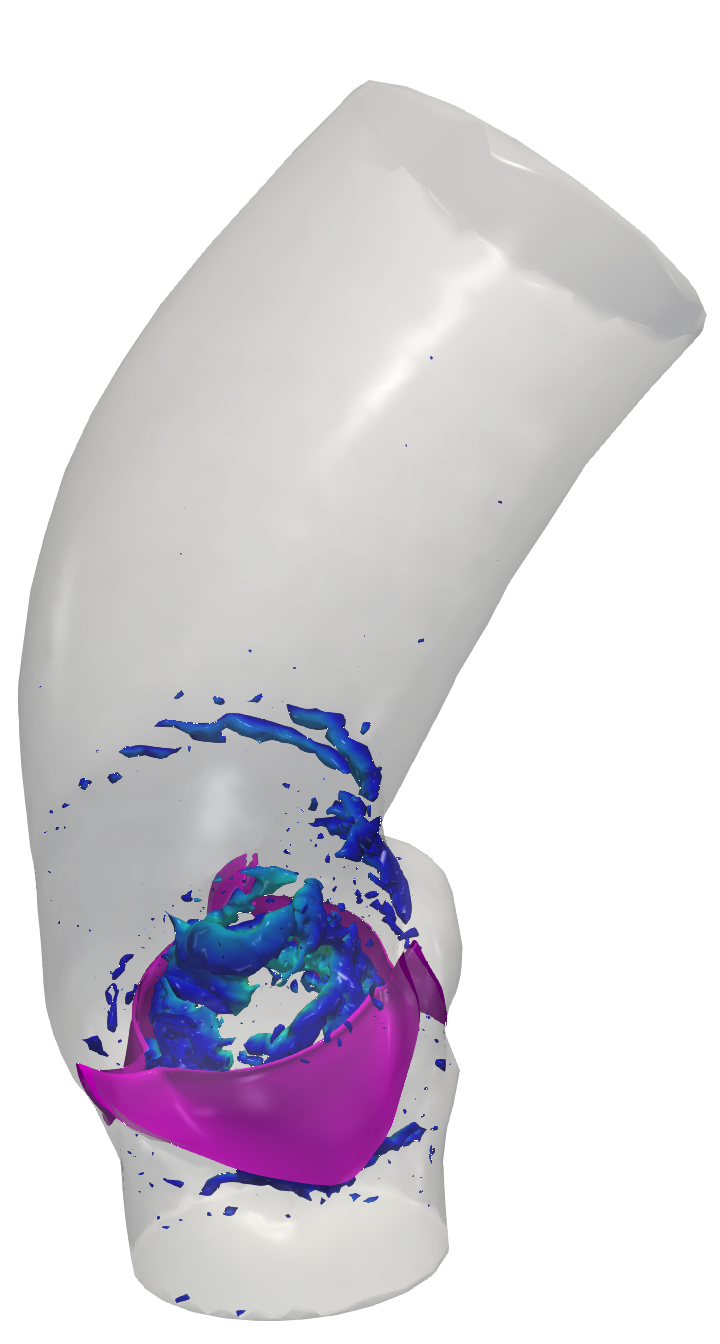}
&
\includegraphics[width=0.17\textwidth]{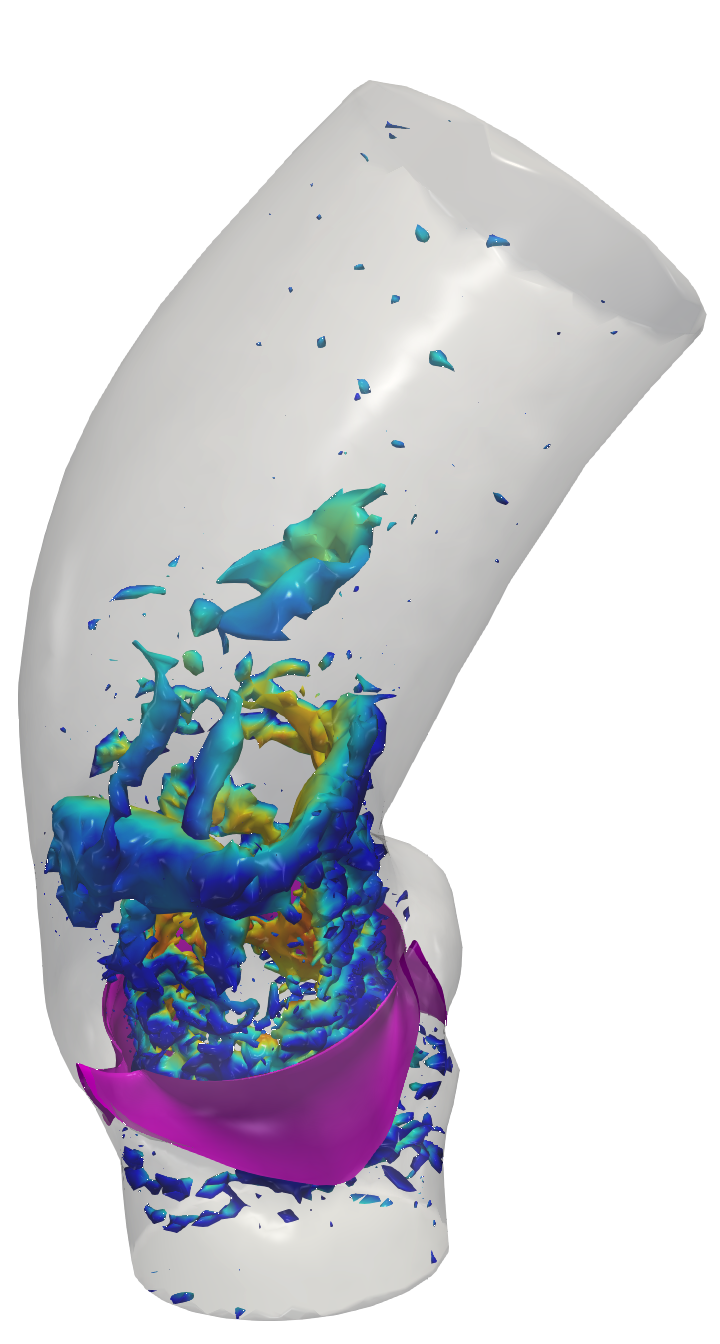}
&
\includegraphics[width=0.17\textwidth]{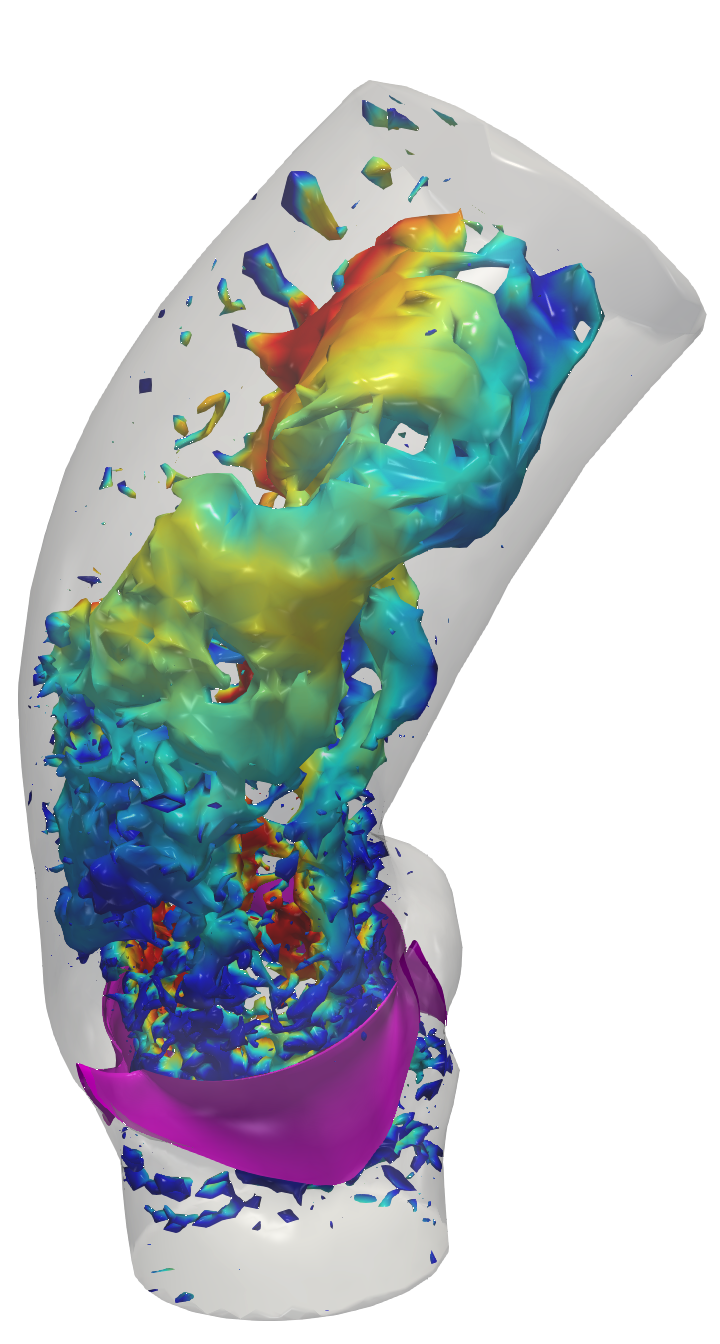}
&
\includegraphics[width=0.17\textwidth]{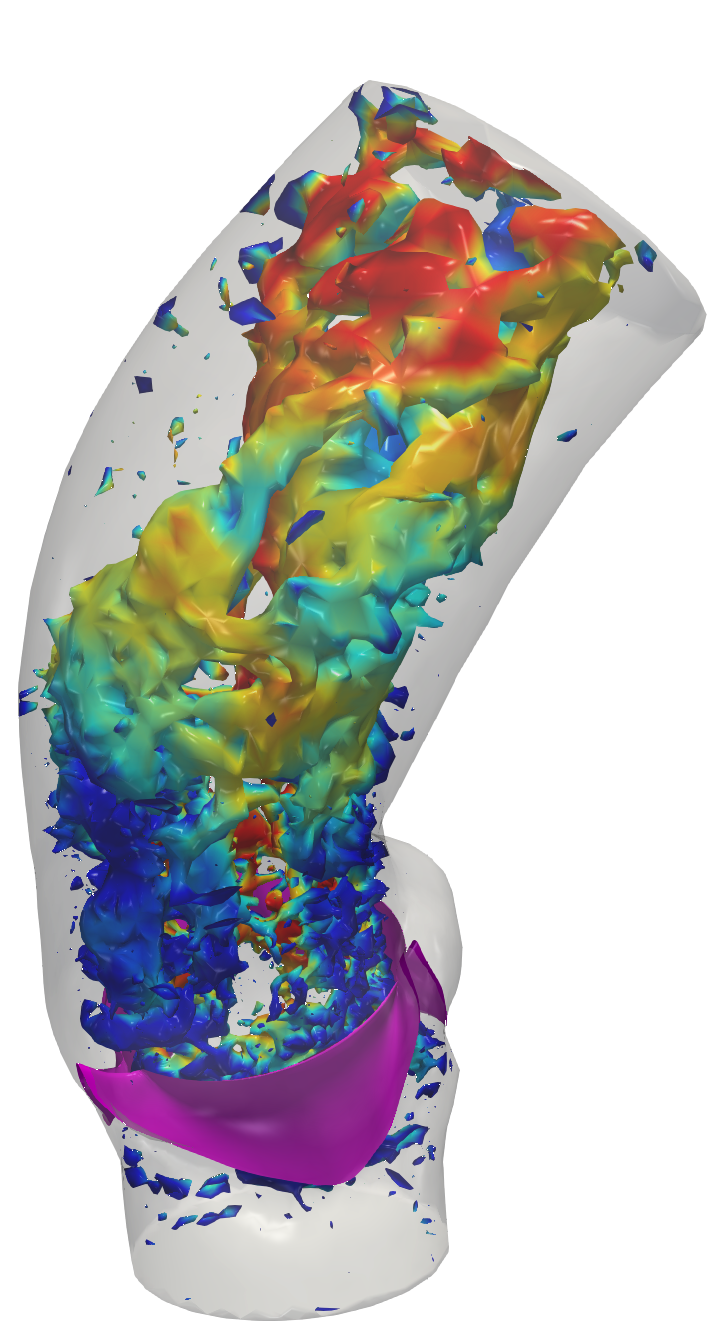}
&
\includegraphics[width=0.17\textwidth]{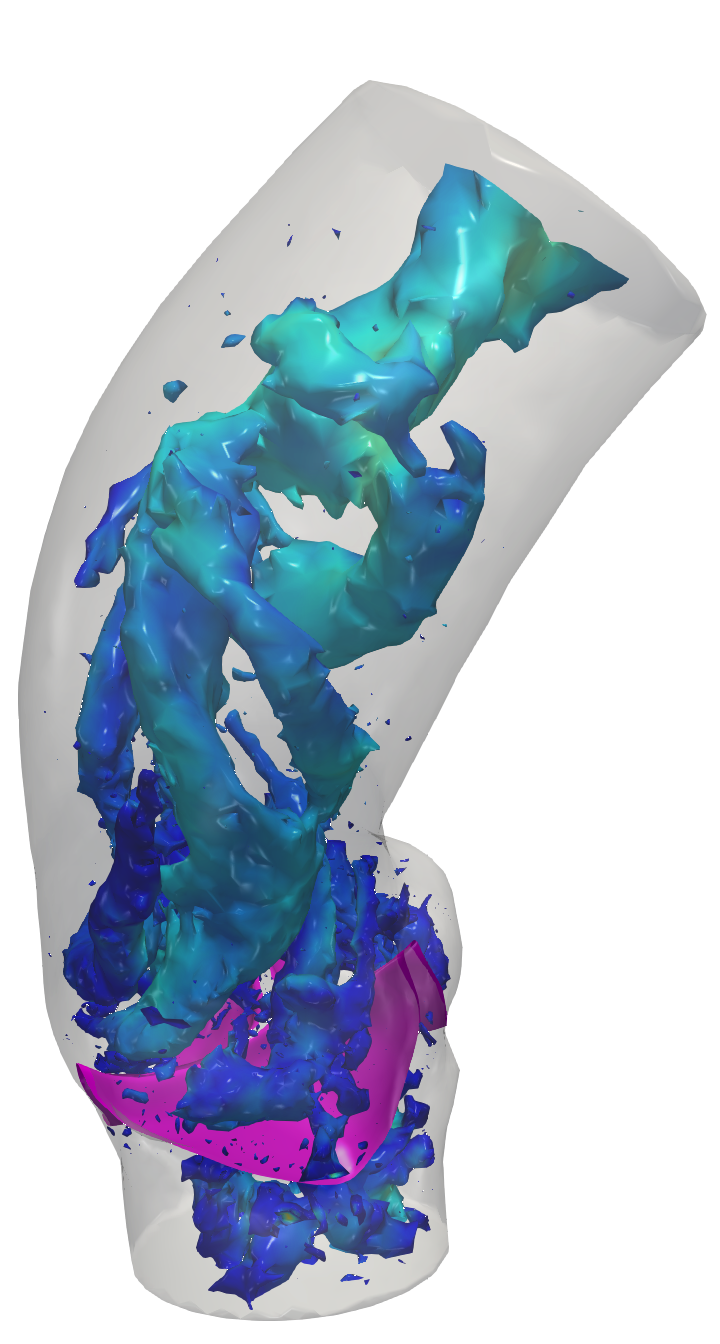}
\\
b) &
\includegraphics[width=0.17\textwidth]{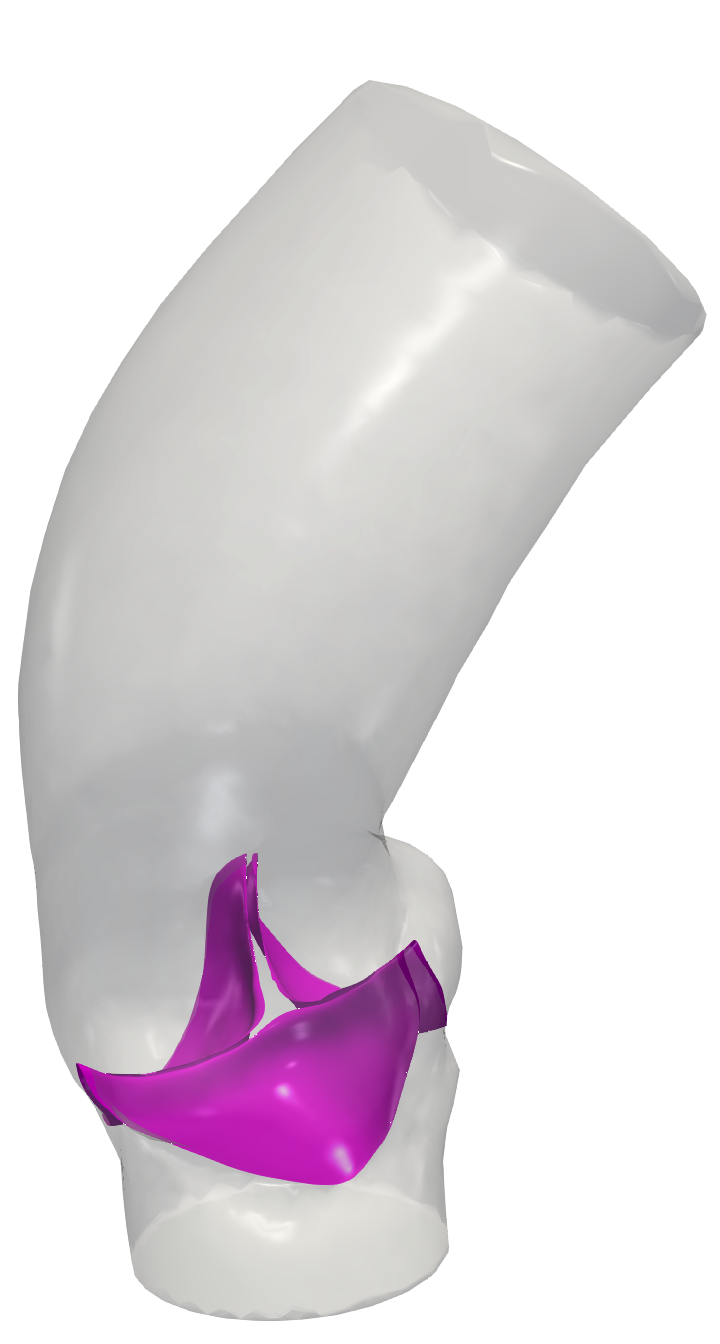}
&
\includegraphics[width=0.17\textwidth]{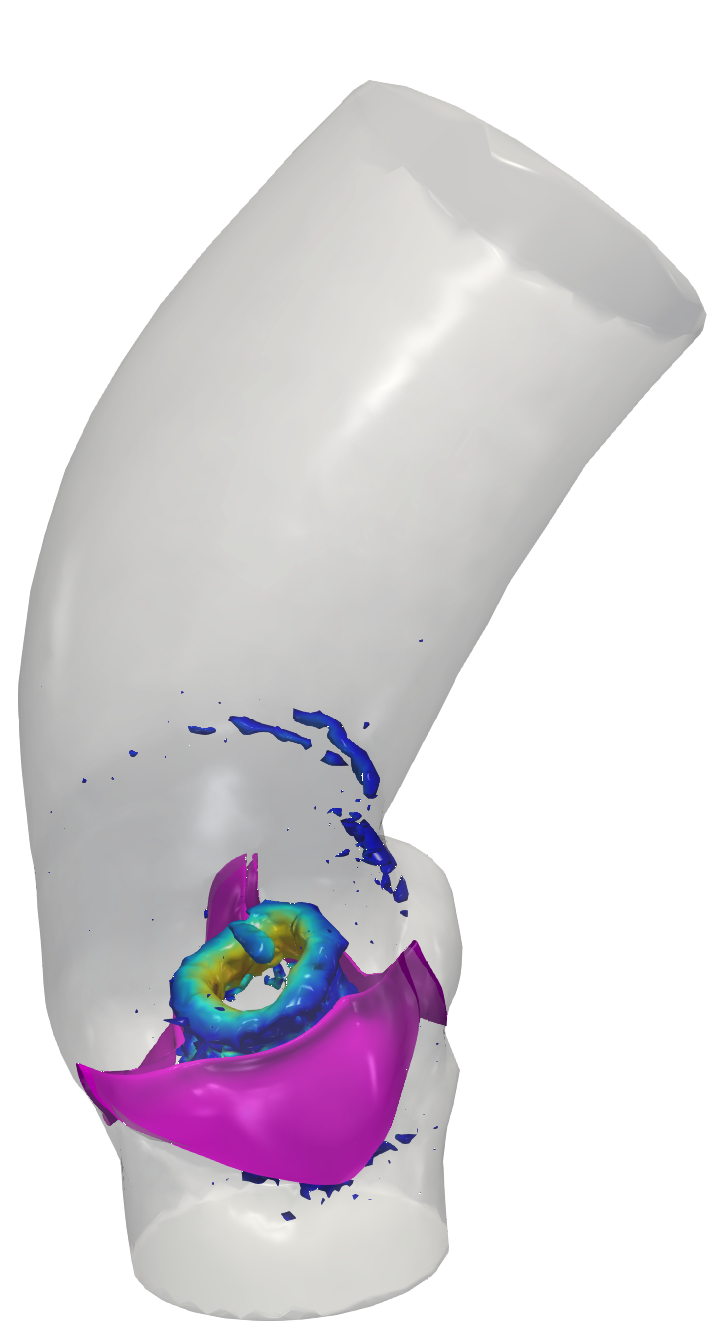}
&
\includegraphics[width=0.17\textwidth]{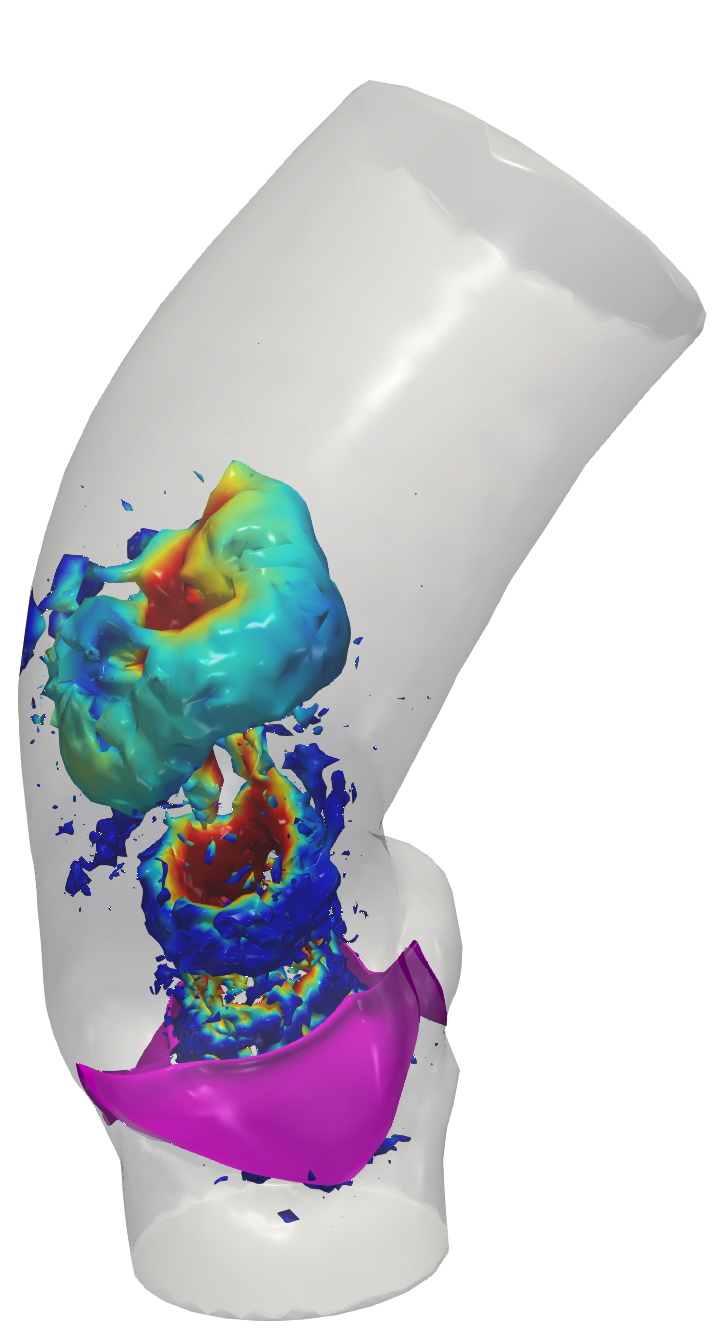}
&
\includegraphics[width=0.17\textwidth]{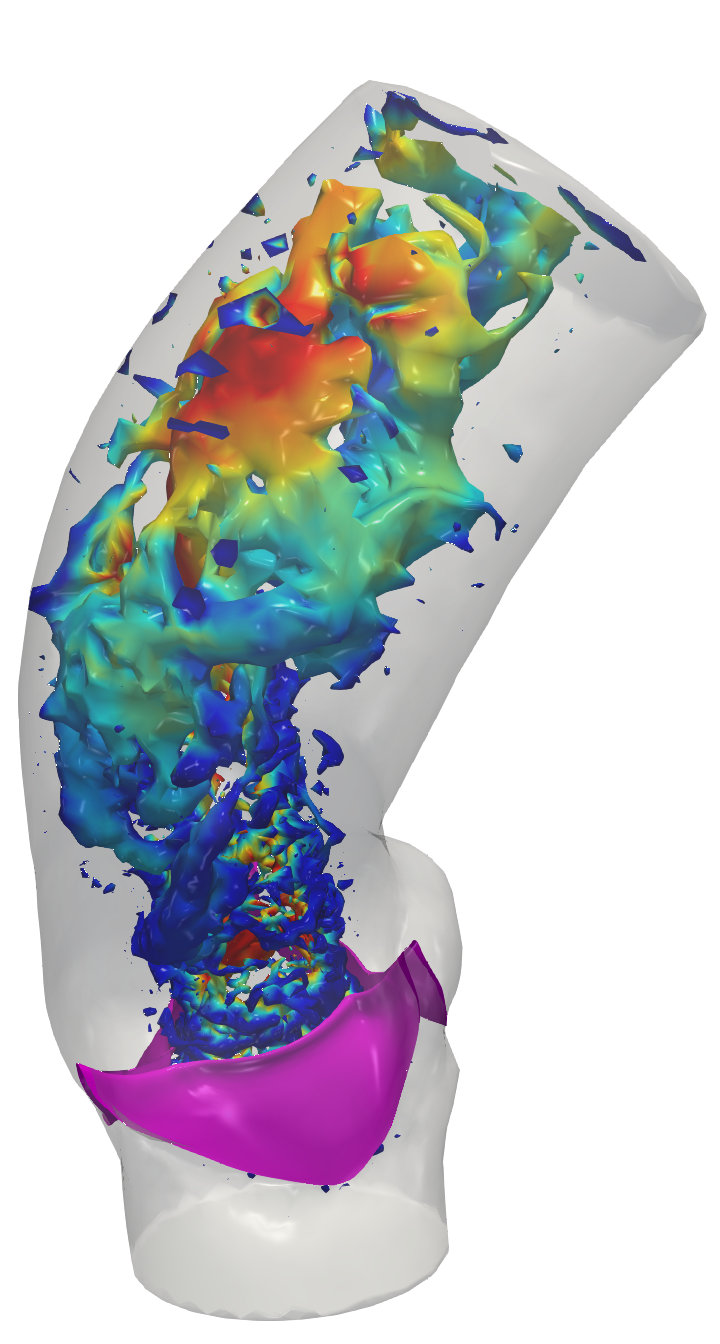}
&
\includegraphics[width=0.17\textwidth]{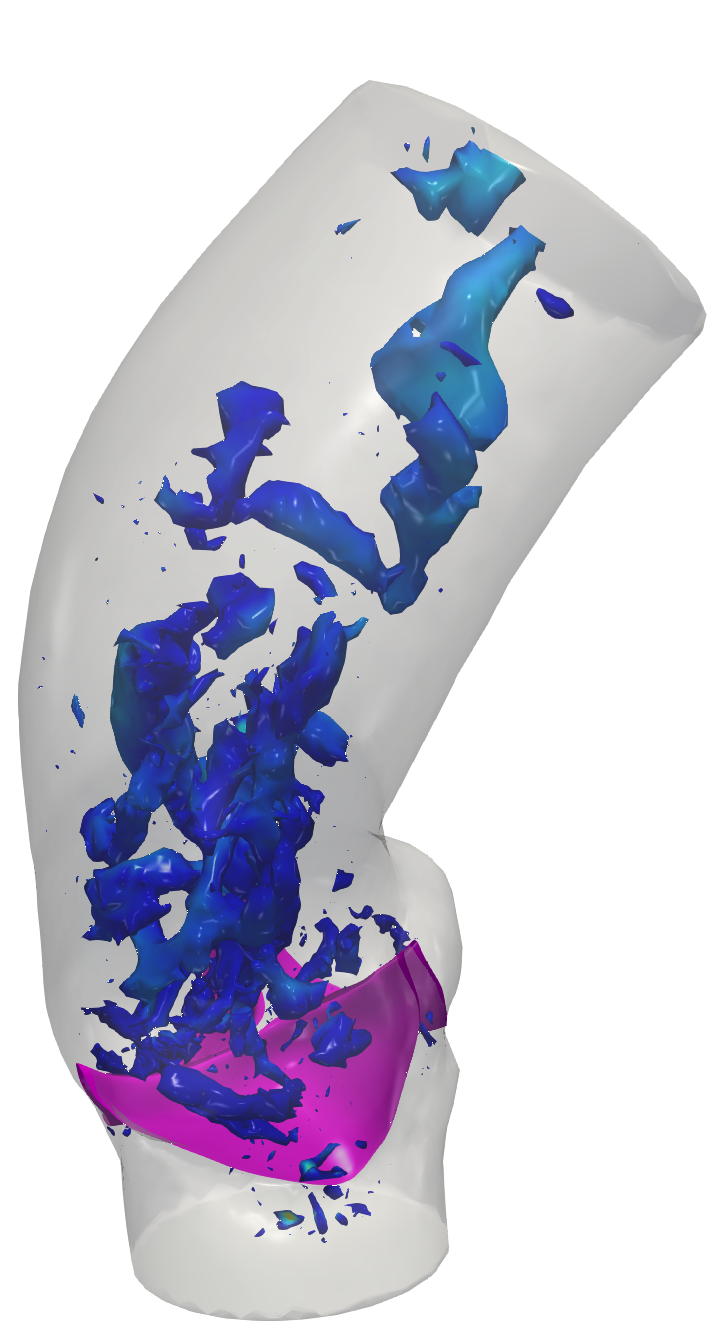}
\\
c) &
\includegraphics[width=0.17\textwidth]{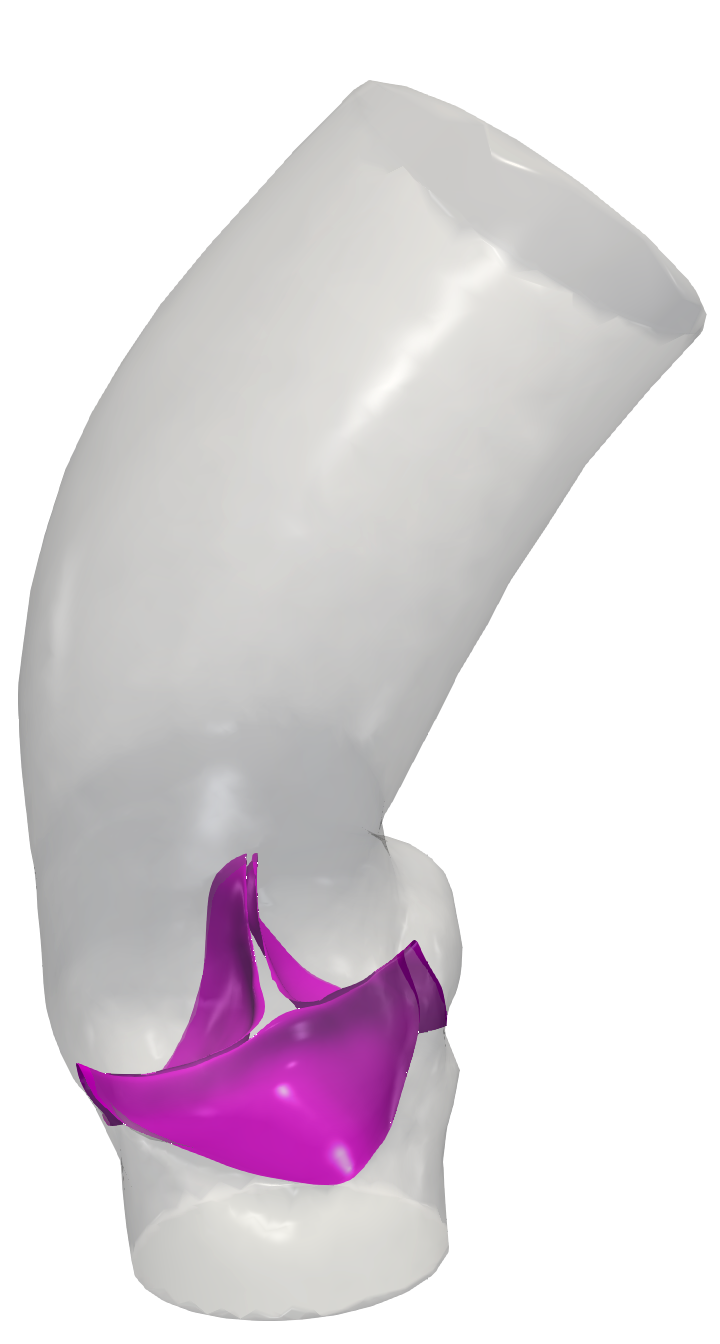}
&
\includegraphics[width=0.17\textwidth]{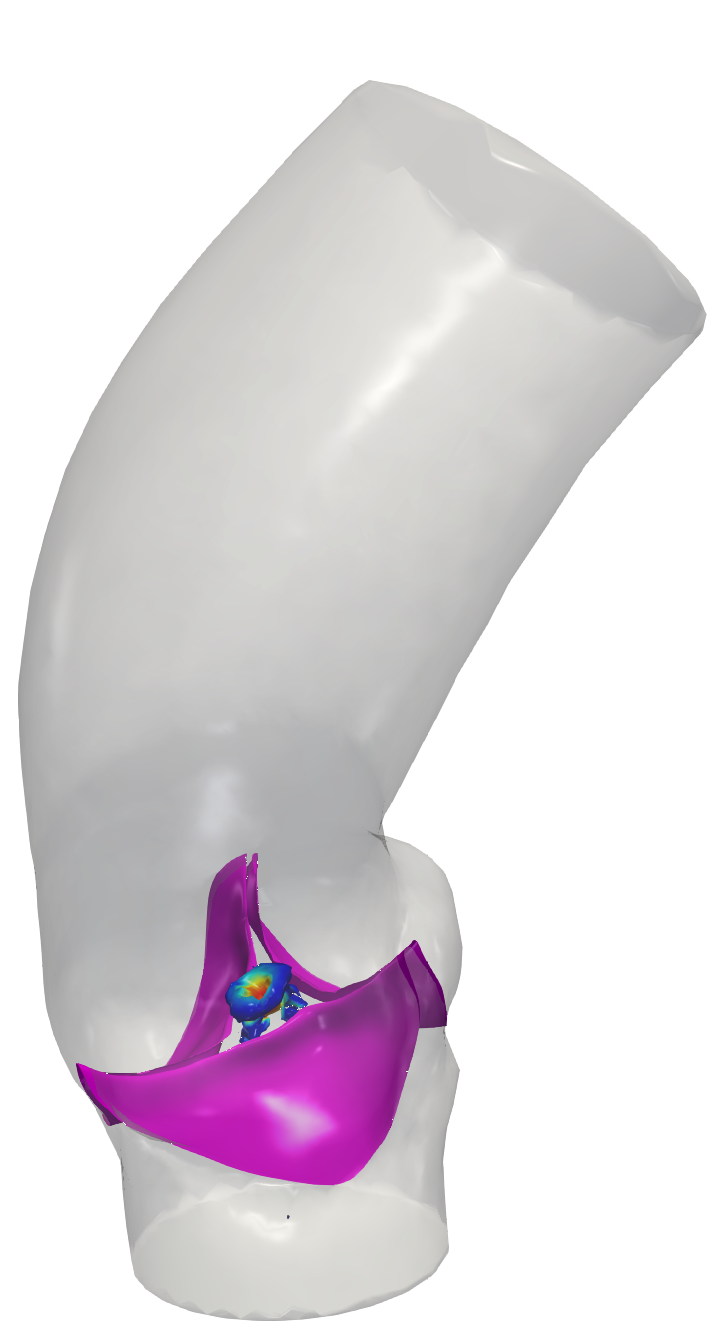}
&
\includegraphics[width=0.17\textwidth]{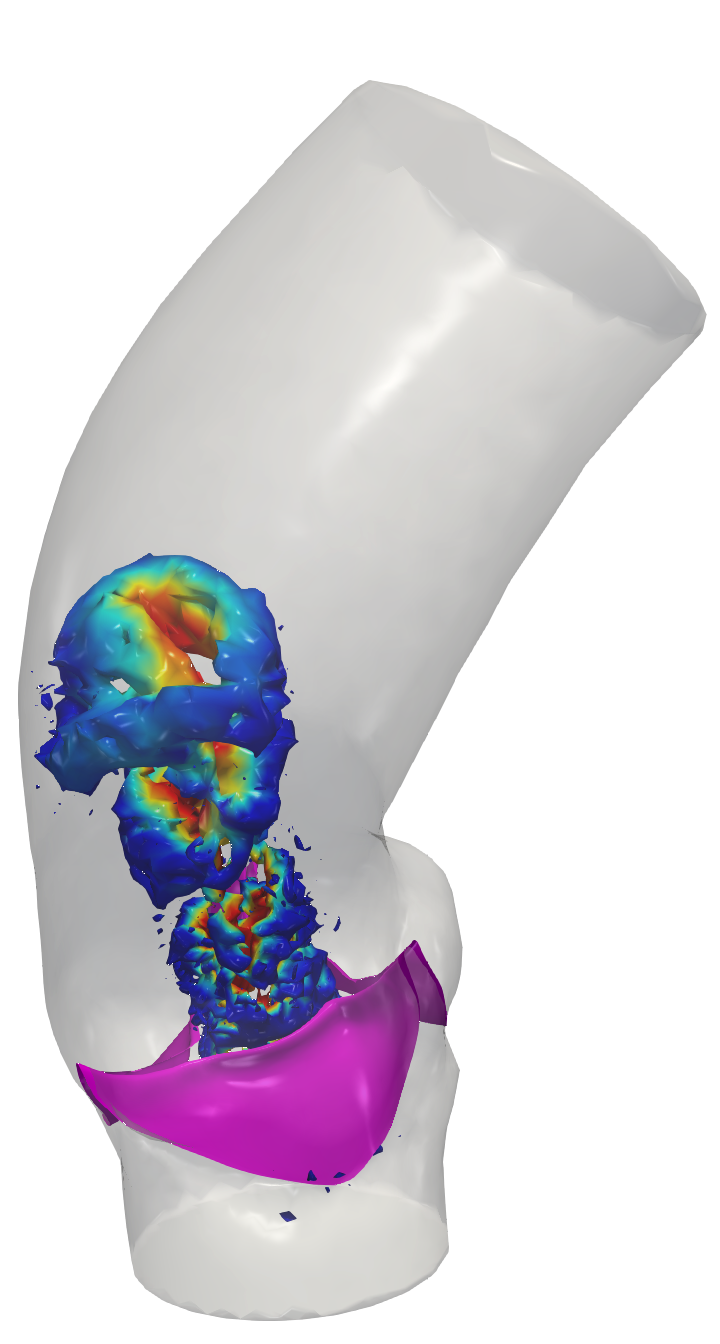}
&
\includegraphics[width=0.17\textwidth]{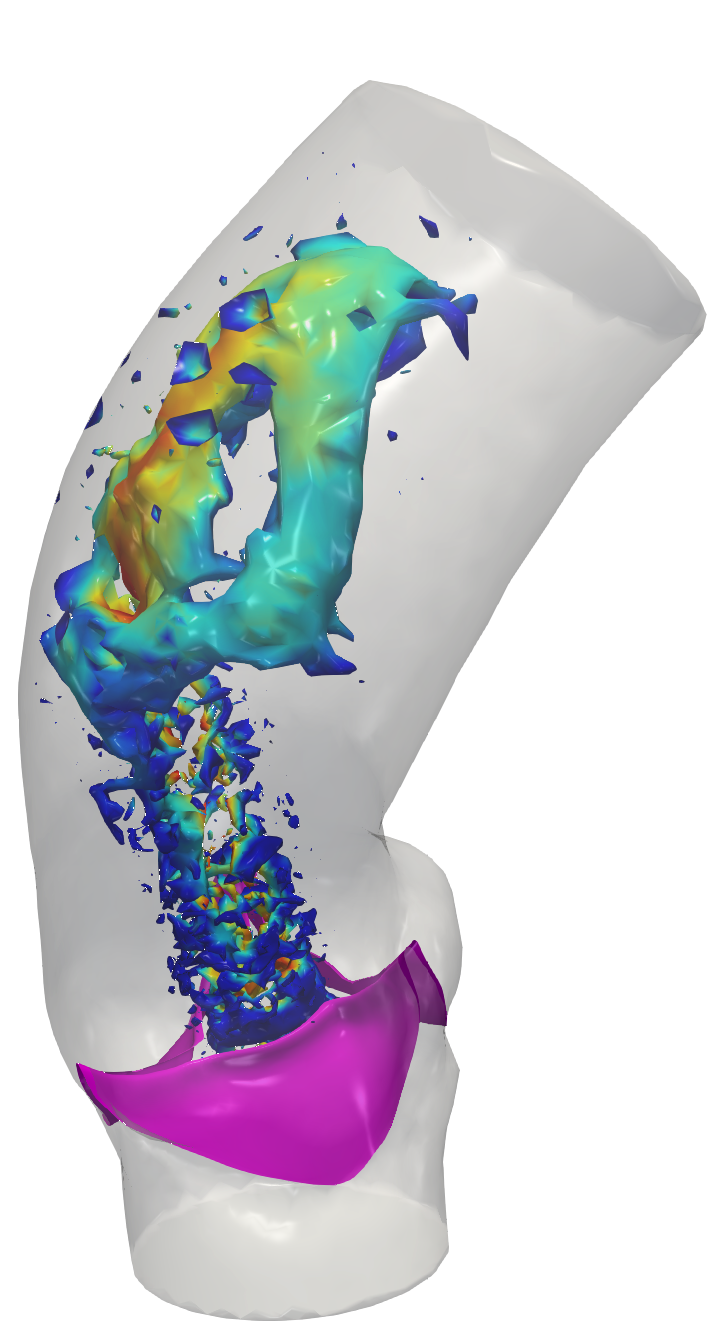}
&
\includegraphics[width=0.17\textwidth]{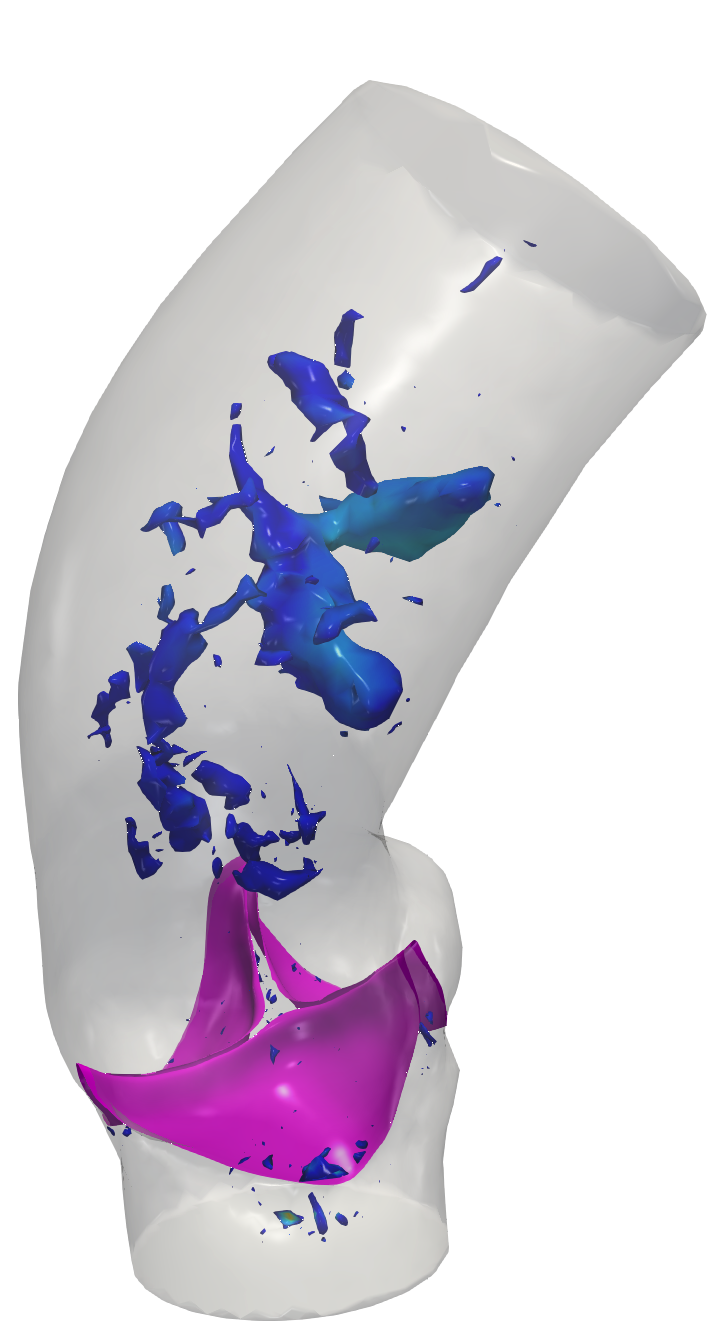}
\\
&{$t=0.15$ s} & {$t=0.20$ s} & {$t=0.25$ s} & {$t=0.30$ s} & {$t=0.40$ s} \\
\end{tabular}\\[1ex]
\ \qquad\qquad\includegraphics[width=0.4\textwidth]{uLegend0to2-hor.png}\caption{
Q-criterion isosurfaces with $Q={5000}$ s\textsuperscript{-2} colored with velocity magnitude at different times: a) physiological case, b) {case \emph{steno-1}}, c) {case \emph{steno-2} (see \cref{tab:synth})}.}
\label{fig:Qcrit}
\end{figure*}

The velocity distibution and the associated coherent vortex structures at different times are displayed in \cref{fig:uSLIC,fig:Qcrit}, respectively.
In the valve opening phase, a jet flow is generated, which leads to the formation of the classical ring coherent structures detaching from the tips of the aortic leaflet (see, e.g., \cite{ringvortex,sotiropoulos2016fluid,becsek2020turbulent}), as we can see at $t=0.2$ s in \cref{fig:Qcrit}.
The vortex structures are then transported downwind in the ascending aorta {during the valve opening phase} and the jet breaks up 
{as soon as the valve is fully open (see \cref{fig:Qcrit}, $t=0.2-0.3$ s)}.
Finally, after the valve is closed, residual flow recirculations can be appreciated both upstream and downstream to the valve.

{
To assess the effectiveness of the RIIS penalty method in representing a non-leaking valve, in \cref{fig:QAVfullSystole} we report the flowrate $Q_\text{AV}$ through a transversal section of the whole domain, together with a zoom on $p_\text{jump}$ when the valve is closed.
By comparison with \cref{fig:cFullSystole}, we can notice that the flowrate is very small when the valve is closed: the maximum of $|Q_\text{AV}|$ when $c=0$ corresponds to a spurious regurgitation of $4.2$ ml/s (attained at $t=0.41$ s), when the valve sustains a negative pressure jump $p_\text{jump}|_{t=0.41\text{s}} \simeq {-98}$ mmHg{, comparable with \cite{hsu2015dynamic}}.
Before the valve is fully closed, instead, in the last part of the closing phase, we can observe a backflow that reaches 296 ml/s: this is due to the valve inertia and it is in partial accordance with the backflow of $\sim200$ ml/s observed in the same phase in \cite{hsu2014fluid,bazilevsHughesFd}, where a detailed 3D valve model is considered.
}
{
\begin{figure}
\centering
\includegraphics[height=0.31\textwidth]{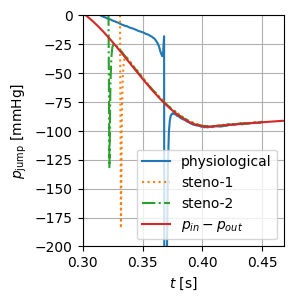}
\includegraphics[height=0.31\textwidth]{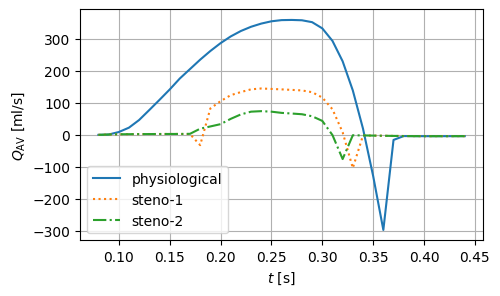}
\includegraphics[width=0.15\textwidth, trim = 0cm -6cm 0cm 0cm]{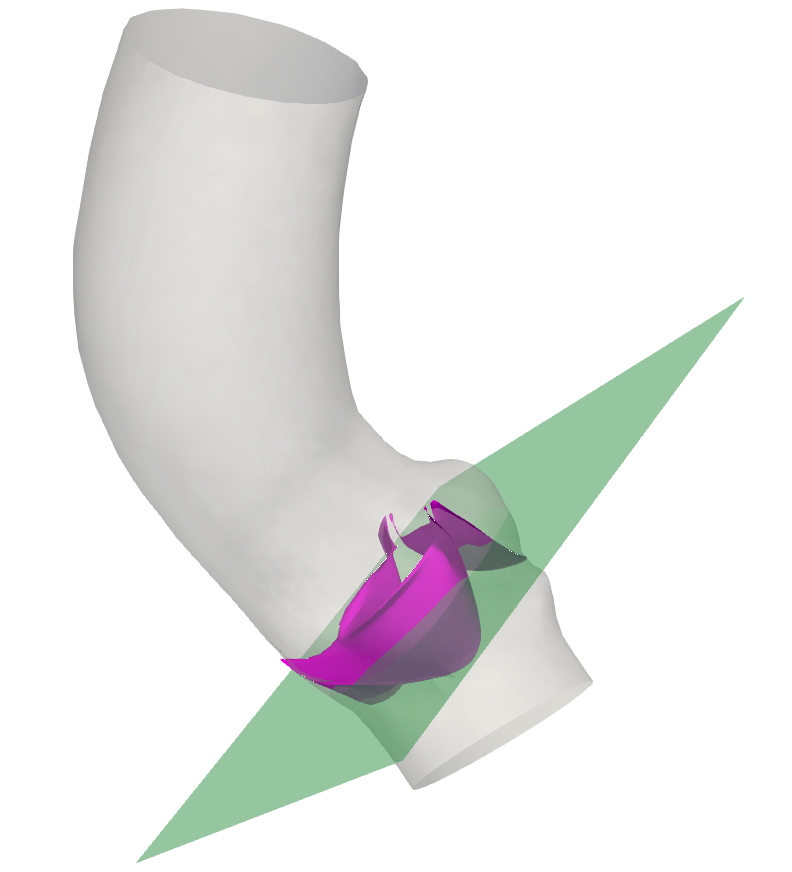}
\caption{{Evolution of $p_\text{jump}$ when the valve is closing and fully closed (left), and of the flowrate $Q_\text{AV}$ throughout the whole systole (center). On the right, the section through which $Q_\text{AV}$ is computed, transversally crossing the whole domain.}}
\label{fig:QAVfullSystole}
\end{figure}
}

\subsubsection{Modeling a stenotic valve}
According to the literature, a calcification-based stenosis of the aortic valve is associated to a reduced compliance of the leaflets, which thus oppose a higher resistance to the blood flow (cf., e.g., \cite{AS}).
This feature can be included in our reduced model by increasing the stiffness parameter $\gamma$, than can be used to model different degrees of stenosis severity.
{Moreover, the added mass of the calcifications increases the leaflets' inertia.}
In this section, we compare the physiological valve with the {the cases \emph{steno-1} and \emph{steno-2} introduced at the beginning of \cref{sec:fullsystole}}. The chosen values of $\gamma$ {and $\rho_\Gamma$} are reported in \cref{tab:synth}
together with the synthetic indicators defined above.

The time evolution of the opening coefficient $c$ and of the orifice area, displayed in \cref{fig:cFullSystole} and summarized in \cref{tab:synth}, shows how {both} an increase in $\gamma$ {and in $\rho_\Gamma$} yields a reduction of the maximum achievable opening of the valve, with an opening phase that is slower {in case \emph{steno-1}.}
Indeed, from \cref{fig:pQFullSystole} we can observe that $p_\text{jump}$ {in the stenotic cases has an almost doubled average value than the physiological cases, and it} even exceeds, in the closing phase, the {maximum} pressure difference $p_\text{in}-p_\text{out}$ that is imposed as a boundary condition over the whole domain.
As already pointed out in \cref{sec:resPhysio}, the control-volume-based, macroscopic pressure jump $p_\text{jump}$ stems from a combination of the \emph{local} stress term $\int_\Gamma\mathbf f\cdot\normal_\Gamma$
and the elastic term $\gamma\int_\Gamma(\curv-\widehat\curv)$.
Indeed, the stress term is comparatively small in the closing phase, which is dominated by the elastic forces, whilst in the opening phase it shows great variability among the three different settings considered here.

In terms of peak velocity $U_\text{peak}$, we notice higher values in the stenotic cases, {with values of more than 2 m/s and thus exceeding the physiological range.
Yet, we point out that,} 
as confirmed by the routine clinical practice (see, e.g., \cite{OA}), a single indicator for stenosis may not be sufficient to fully categorize a patient's condition, and different indicators have to be considered at the same time.
{Indeed, although the values of AS may indicate that {steno-1} and {steno-2} are representative of mild and moderate stenosis, respectively (see \cite{mittal0d,OA}), the indicators $U_\text{peak}$ and $p_\text{jump,peak}$ both lie in the mild stenosis range.}

In \cref{fig:uSLIC,fig:Qcrit} we report the velocity field on a 2D slice and the coherent vortex structures generated by the Q-criterion method (cf.~\cite{Qcriterion}).
We notice that, in the stenosis cases, the reduced orifice area and the shorter time interval in which flow is allowed through the valve yield a stronger aortic jet and a more disorganized velocity distribution.
The ring vortex detaching from the tips of the leaflets, that can be seen at time $t=0.20$ s, is highly distorted in a short time (see {$t=0.25-0.30$} s) while {it is transported along the jet and impacts} on the posterior aortic wall.
{After valve closure (see $t=0.40$ s), the smaller blood velocity magnitude and vortical structure dimensions in the stenotic cases indicate a less effective mixing of blood and a longer residence time of blood in this portion of the vessel, thus yielding a reduced cardiac output.}
{Moreover, since \emph{steno-2} represents a more stenotic case than \emph{steno-1}, }
the jet that can be appreciated in \cref{fig:uSLIC}, c), albeit characterized by high velocity values, gets very narrow and lasts for {less than}
half of the systole.
Correspondingly, 
{the velocity profile is more chaotic, and the vortical structures undergo a faster breakdown}
into small-scale eddies {(cfr.~\cref{fig:Qcrit}, b)-c)).}

\subsubsection{Comparison with the Korakianitis-Shi model}\label{sec:KS}
We compare the proposed curvature-based reduced model with another 3D-0D FSI system proposed by \cite{resistivo}.
In that reference, a Navier-Stokes-RIIS fluid dynamics system is coupled with the Korakianitis-Shi (KS) model for {leaflet} mechanics introduced by \cite{korshi}:
\begin{equation}
\ddot\theta + k_f\dot\theta = (k_p p_\text{jump}+k_bQ)\cos\theta - k_v\,|Q|\,\sin(2\theta),
\end{equation}
where $\theta$ is the valve opening angle, ranging between some prescribed values $\theta_\text{min},\theta_\text{max}$, and it depends on the pressure jump $p_\text{jump}$ and the flowrate $Q$ across the valve; $k_{(\cdot)}$ are model parameters that need calibration.
Hinging upon the definition of the resistance area introduced by \cite{korshi,resistivo}
\begin{equation}
\text{AR}_\text{ao} = \frac{(1-\cos\theta)^2}{(1-\cos\theta_\text{max})^2}
\end{equation}
and observing that the orifice area is quadratic {with respect to} the opening coefficient $c$ ranging from 0 to 1, the opening angle $\theta(t)$ can be related to $c(t)$ by
\begin{equation}
c(t) = \frac{\cos\theta_\text{min} - \cos\theta(t)}{\cos\theta_\text{min} - \cos\theta_\text{max}}.
\end{equation}

\begin{table}
\centering
\begin{tabular}{cccccc}
$k_f$ & $k_p$ & $k_b$ & $k_v$ & $\theta_\text{min}$ & $\theta_\text{max}$\\
$\left[\frac{\text{rad}}{\text{s}}\right]$ & $\left[\frac{\text{rad}}{\text{mmHg s}^2}\right]$ & $\left[\frac{\text{rad}}{\text{mL}}\right]$ & $\left[\frac{\text{rad}}{\text{mL}}\right]$ & $[\,^\circ\,]$ & $[\,^\circ\,]$\\
\hline
50 &
{665}
& 6 & 21 & 5 & 75
\end{tabular}
\caption{Model parameters for the KS model.\label{tab:KS}}
\end{table}

In reference \cite{resistivo}, a patient-specific geometry was analyzed, and a quasi-static approach was adopted in the RIIS term, that is $\mathbf u_\Gamma=\mathbf 0$ is considered in the momentum equation.
In view of the discussion of \cref{sec:quasistatic}, we drop here the quasi-static hypothesis: this and the difference in the geometry of interest lead us to a re-calibration of the KS model.
Consistently with what done for the proposed curvature-based 0D model, the calibration was carried out aiming at a physiological opening time, and the resulting model parameters are reported in \cref{tab:KS}.

\begin{figure}
\centering
\includegraphics[width=0.47\textwidth]{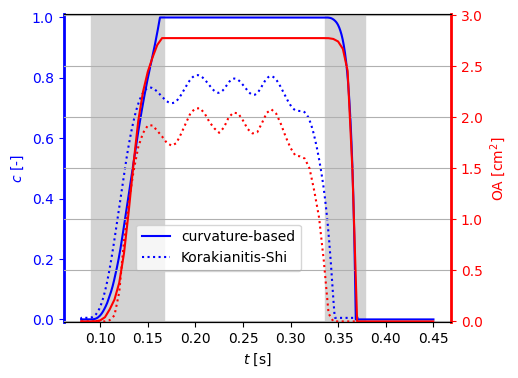}
\caption{{Full systole: opening coefficient $c$ (left axis) and orifice area OA (right axis) under physiological pressure conditions, obtained with the curvature-based model and with the KS model.}}
\label{fig:cFullSystoleKS}
\end{figure}

\begin{table}
\centering
\begin{tabular}{rl|cc}
\multicolumn{2}{c|}{model} & curvature-based & KS
\\
\hline
$\gamma$ & [N/m]
    & {3}
    & -
\\
$T_\text{open}$ & [ms]
    & {72}
    & {64}
\\
$T_\text{close}$ & [ms]
    & {32}
    & {52}
\\
$OA_\text{max}$ & [cm\textsuperscript{2}]
    & {2.78}
    & {2.10}
\\
$AS$ & [\%]
    & 0
    & {24}
\\
$U_\text{peak}$ & [m/s]
    & {1.52}
    & {1.78}
\\
$p_\text{jump, peak}$ & [mmHg] & {3.45} & {5.98}
\end{tabular}
\caption{Synthetic indicators for valve stenosis. \label{tab:synthKS}}
\end{table}

The time evolution of the opening coefficient and the associated effective orifice area are displayed in \cref{fig:cFullSystoleKS}, and synthetic indicators are reported in \cref{tab:synthKS}.
We notice that the resulting KS model does not allow a full opening of the valve, with a maximum angle $\theta=66^\circ<\theta_\text{max}$,{although both the opening and closing times $T_\text{open}$ and $T_\text{close}$ lie in the physiological ranges $76\pm 30$ ms and $46\pm 12$ ms, respectively}.
The latter observation may be seen as an improvement with respect to \cite{resistivo}, which reported a slow closing phase{, that can be ascribed to considering the velocity surface $\mathbf{u}_\Gamma$.}
The lack of reaching a fully open position, in turn, is due to the choice made on the calibration strategy: in additional numerical tests we observed that modifying the parameters to attain a larger maximum value of $\theta$ would determine a non-physiologically short opening time (an observation in accordance with \cite{resistivo}).

\begin{figure*}
\centering
\begin{tabular}{m{1ex}ccccc}
a) &
\includegraphics[width=0.17\textwidth]{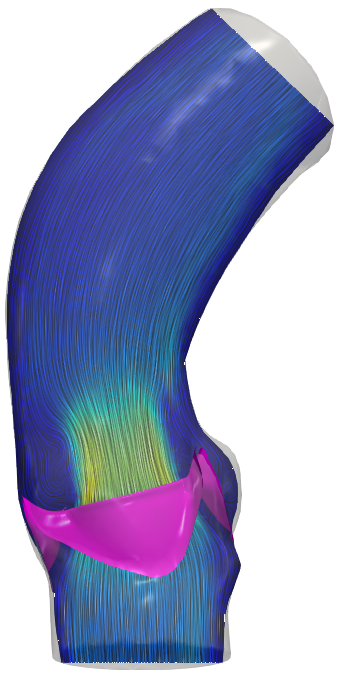}
&
\includegraphics[width=0.17\textwidth]{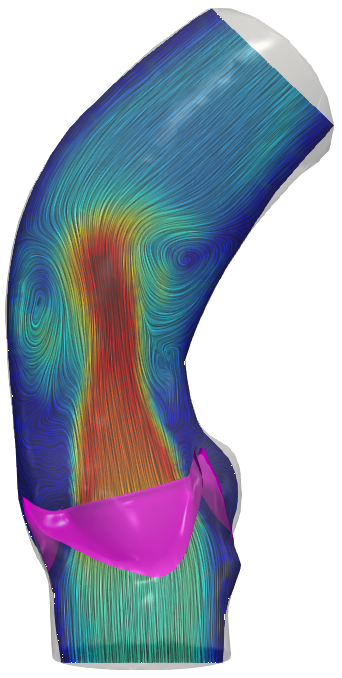}
&
\includegraphics[width=0.17\textwidth]{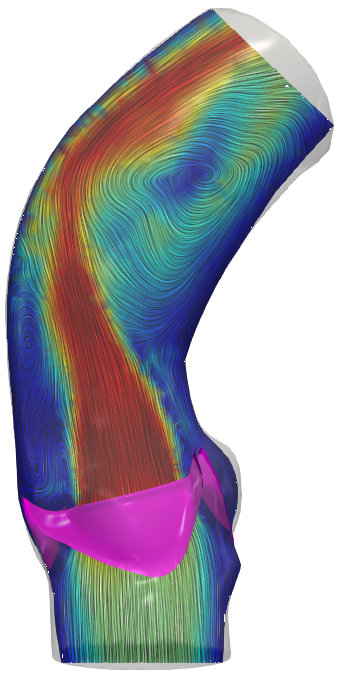}
&
\includegraphics[width=0.17\textwidth]{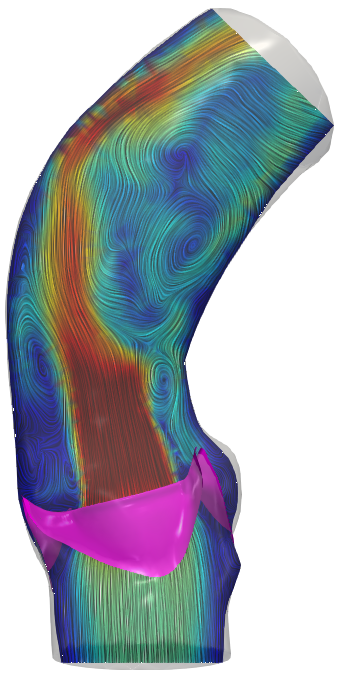}
&
\includegraphics[width=0.17\textwidth]{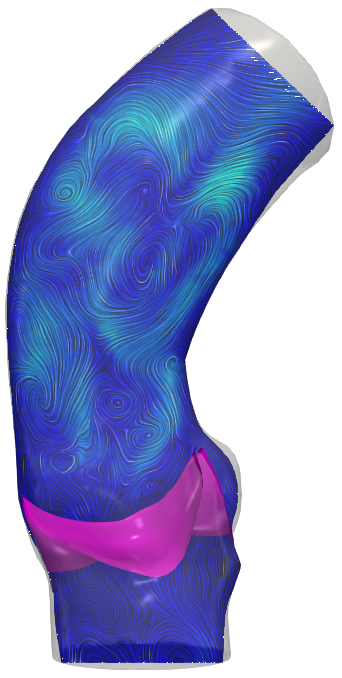}
\\
b) &
\includegraphics[width=0.17\textwidth]{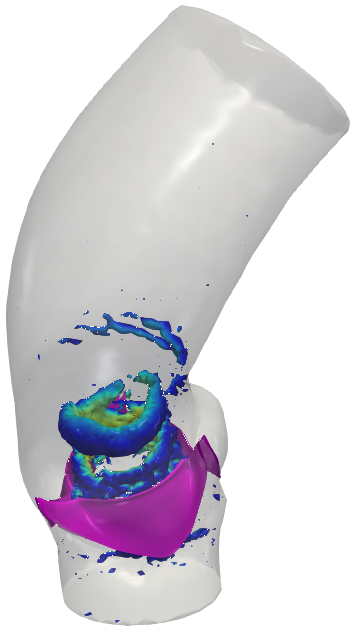}
&
\includegraphics[width=0.17\textwidth]{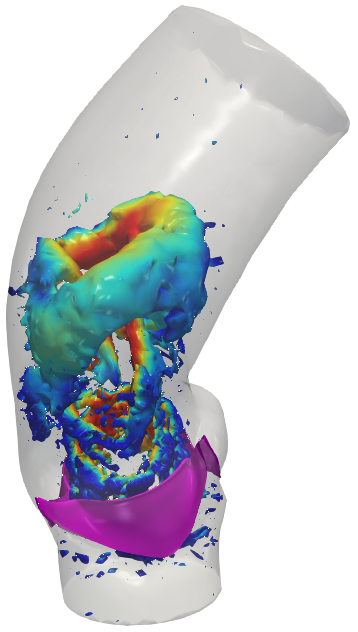}
&
\includegraphics[width=0.17\textwidth]{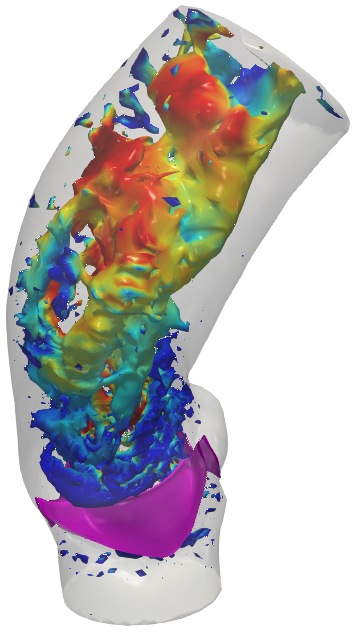}
&
\includegraphics[width=0.17\textwidth]{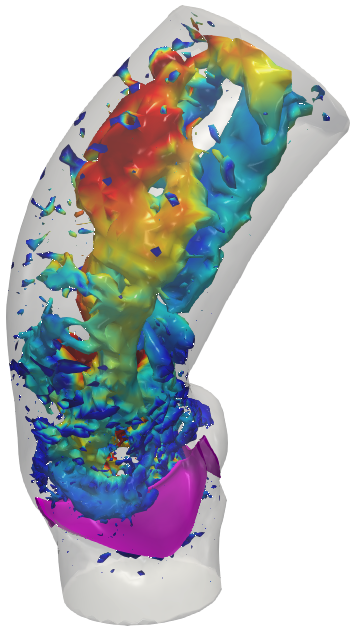}
&
\includegraphics[width=0.17\textwidth]{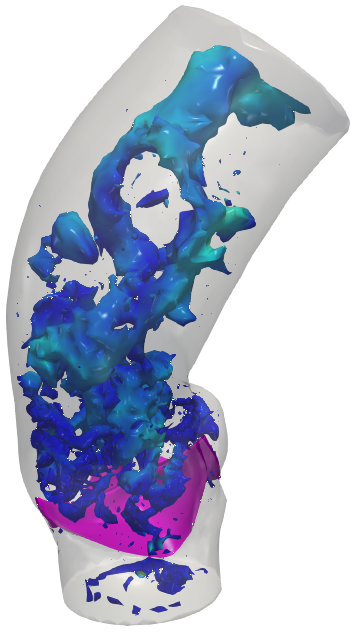}
\\
&$t=0.15$ s & $t=0.20$ s & $t=0.25$ s & $t=0.30$ s & $t=0.40$ s \\
\end{tabular}\\[1ex]
\ \qquad\qquad\includegraphics[width=0.45\textwidth]{uLegend0to2-hor.png}
\caption{
Results with Korakianitis-Shi model at different times:
a) velocity distribution on a longitudinal slice;
b) Q-criterion isosurfaces with $Q={5000}$ s\textsuperscript{-2} colored with velocity magnitude.}
\label{fig:uSLICandQcritKS}
\end{figure*}

In \cref{fig:uSLICandQcritKS} we report the velocity distribution and coherent vortical structures obtained in these settings.
Comparing such results with those of \cref{fig:uSLIC,fig:Qcrit},
{we notice that that they are intermediate between the physiological and \emph{steno-1} cases of the curvature-based model, as it is for the stenosis indicators $AS$ and $U_\text{peak}$ reported in \cref{tab:synthKS}.}
This confirms that the KS model, when calibrated in order to attain physiological opening times, leads to a slightly stenotic behavior of the valve.

\section{Conclusions}

We proposed a novel reduced FSI model for the aortic valve.
The valve dynamics was described by a lumped-parameter model considering the flow-induced stress and a curvature-based elasticity term, as well as damping effects, and its coupling with the 3D blood flow was based on the RIIS method.
This system was employed to simulate the blood flow in the ascending aorta, both in physiological conditions and in the case of {mild} aortic stenosis.

The numerical results {demonstrate that the proposed model is a computationally efficient approach for simulating aortic hemodynamics and the effects of valve dynamics on blood flow.
Compared to a CFD simulation with prescribed leaflet displacement, the additional computational effort cost is minimal, limited to the assembly of the right-hand side of the ODE governing valve dynamics, which can be efficiently carried out at each quadrature node.}
The model also straightforwardly provides an explicit expression for the leaflets' velocity $\uu_\Gamma$, without resorting to complex reconstruction procedures that would introduce discrete interpolation errors.
The comparison with a quasi-static approach adopted in previous works and with the Korakianitis-Shi model showed the advantages of our model in reconstructing the surface velocity and reproducing a physiological duration of the valve opening and closing phases.

{
To achieve the computational efficiency of the proposed method, we introduced some assumptions that may limit its applicability.
Although accounting for macroscopic curvature changes in the leaflets, the model does not fully describe local deformations nor leaflet coaptation or prolapse.
Moreover, we assume uniform material properties, ignoring the heterogeneity of stiffness and thickness and the anisotropy of the tissue, which may impact valve dynamics.
For these reasons, in several applications, a reduced model such as the one proposed is not a substitute for a fully three-dimensional fluid-structure model.
The latter approach is necessary to have an accurate description of mechanical stresses in the leaflets and a locally detailed stress-strain relationship, as needed, e.g., in the investigation of the onset and progression of valve calcification or structural degeneration.
Moreover, a fully detailed modeling is required to analyze flow details in proximity of the leaflets and shear stress distributions, associated with thrombotic risk,
or to capture leaflet fluttering, which is the subject of increasing investigation especially in prosthetic valve design. 
Finally, our reduced model does not allow to predict long-term biomechanical changes related to valve disease progression, remodeling or degeneration.}

{Nevertheless, our model could be employed as an agile computational tool for several hemodynamics investigations. 
Since it includes the valve geometry, it can provide a more realistic representation of transvalvular pressure gradients and flow features compared to lumped parameter models.
It could be used in scenario analysis for the assessment of the hemodynamics effects of valve stenosis in domains including the ventricle or a wider downstream tract of the aorta and proximal vessels, where the computational burden of a fully 3D model may become impractical.
Indeed, a preliminary step in this direction has been taken in the investigation of pulmonary valve replacement \cite{criseo2024computational}.
Furthermore, it could be employed in population-based studies, where computational efficiency is paramount.
}

In order to further enhance the proposed model for the study of different scenarios and pathological conditions, different directions of research may be undertaken.
{
A more precise representation of the valve's open configuration could be achieved using patient-specific imaging data. Moreover, the opening field $\dbar$ could be replaced by a more complex displacement field involving additional (though still limited) degrees of freedom, as in \cite{pedrizzetti0d}. This could improve the conservation of the valve's mass throughout its dynamics and possibly account for additional kinematic modes in a synthetic way.
}

An efficient semi-automatic calibration strategy of the model parameters would allow a patient-specific analysis of different pathological conditions, as well as the simulation of possible treatment scenarios to help the pre-operative design.
In such context, the efficiency of the calibration procedure would have particular relevance: reduced order models and machine-learning-based surrogates of this complex system may thus help in this respect.
Thanks to the general derivation of the model from a local force balance, its extension to the pulmonary valve of even the atrio-ventricular valves can be envisaged, possibly introducing additional terms accounting for the subvalvular apparatus.

{In terms of the numerical scheme, an implicit coupling of the fluid and valve model could be considered.
While this would increase computational cost, adopting a semi-implicit strategy, as proposed in \cite{hsu2015dynamic,johnson2020thinner}, could help mitigate the additional effort.}

Finally, an additional level of complexity may be introduced by considering contact forces exchanged among the leaflets, that may affect the dynamics in the early opening phase and in diastole.

\section{Acknowledgments}
LD acknowledges funding from the {HORIZON-EUROHPC-JU-2023-COE-03} project {\bf dealii--X} ``an Exascale Framework for Digital Twins of the Human Body'' (no. 101172493).
We gratefully acknowledge the CINECA award under the ISCRA C initiative, for the availability of high performance computing resources and support in the project HP10C794N5 (P.I. I.~Fumagalli).
{All the authors are members of the INdAM group GNCS ``Gruppo Nazionale per il Calcolo Scientifico'' (National Group for Scientific Computing).
The present research is part of the activities of ``Dipartimento di Eccellenza 2023–2027'', MUR, Italy,
Dipartimento di Matematica, Politecnico di Milano.}


\bibliographystyle{myalpha}
\bibliography{valve0Dmodel}

\end{document}